\documentclass{amsart}

\usepackage{amssymb}
\usepackage{amsmath}
\usepackage{amsthm,amsrefs}
\usepackage{graphicx}
\usepackage{tikz-cd}
\usepackage{enumitem}
\usepackage{tensor}

\usepackage{tikz}
\usetikzlibrary{arrows}
\usetikzlibrary{decorations.markings}
\usetikzlibrary{patterns}

\newcommand{\N}{\mathbb{N}}
\newcommand{\NN}{\mathcal{N}}
\newcommand{\F}{\mathcal{F}}

\newcommand{\PP}{\mathbb{P}}

\newcommand{\Z}{\mathbb{Z}}

\newcommand{\BB}{\mathcal{B}}
\newcommand{\MM}{\mathcal{M}}
\newcommand{\X}{\mathbb{X}}
\newcommand{\Y}{\mathbb{Y}}
\newcommand{\CC}{\mathcal{C}}
\newcommand{\DD}{\mathcal{D}}
\newcommand{\E}{\mathcal{E}}
\newcommand{\cP}{\mathcal{P}}
\newcommand{\cS}{\mathcal{S}}
\newcommand{\A}{\mathcal{A}}

\newcommand{\TT}{\mathcal{T}}

\newcommand{\m}{\mathfrak{m}}

\renewcommand{\L}{\mathcal{L}}
\renewcommand{\H}{\mathcal{H}}
\renewcommand{\O}{\mathcal{O}}

\renewcommand{\mod}{\mathsf{mod}}
\newcommand{\Gr}{\mathsf{Gr}}
\newcommand{\gr}{\mathsf{gr}}
\newcommand{\Hom}{\mathsf{Hom}}
\newcommand{\Ext}{\mathsf{Ext}}
\newcommand{\QGr}{\mathsf{QGr}}
\newcommand{\qgr}{\mathsf{qgr}}
\newcommand{\QHom}{\mathcal{H}om}
\newcommand{\HOM}{\underline{\Hom}}
\newcommand{\QHOM}{\underline{\QHom}}
\newcommand{\QGamma}{\mathbf{\Gamma}}
\newcommand{\End}{\mathrm{End}}
\newcommand{\coh}{\mathsf{coh}}
\newcommand{\coker}{\mathrm{coker ~}}
\newcommand{\Ab}{\mathsf{Ab}}
\newcommand{\K}{\mathsf{K}}
\newcommand{\diag}[1]{{#1}_{\Delta}}
\newcommand{\RHom}{\mathsf{RHom}}
\newcommand{\dtensor}{\otimes^L}

\newcommand{\ssf}[1]{\mathsf{#1}}

\newtheorem{theorem}{Theorem}[section]
\newtheorem*{theorem*}{}
\newtheorem{lemma}[theorem]{Lemma}
\newtheorem{propn}[theorem]{Proposition}
\newtheorem{cor}[theorem]{Corollary}

\theoremstyle{definition}
\newtheorem{defn}[theorem]{Definition}
\newtheorem{ex}[theorem]{Example}
\newtheorem{remark}[theorem]{Remark}

\numberwithin{equation}{section}

\begin{document}

\title{Coherent $I$-indexed algebras and noncommutative projective schemes}

\author{Jackson Ryder}

\address{Department of Mathematics and Statistics, University of New South Wales}
\email{jackson.ryder@unsw.edu.au}


\begin{abstract}
    We study coherent $I$-indexed algebras and associated noncommutative projective schemes, where the index set $I$ is a locally finite directed poset. Our main result is a characterisation of such noncommutative projective schemes in terms of sequences of objects with properties analogous to the sequence of tensor powers of an ample line bundle, extending a similar characterisation given by Polishchuk for coherent $\mathbb{Z}$-indexed algebras. 
\end{abstract}

\maketitle

\section{Introduction}
In \cite[Theorem 2.4]{polishchuk_noncommutative_2005}, Polishchuk extended a theorem of Artin and Zhang \cite[Theorem 4.5]{artin_noncommutative_1994} characterising abelian categories equivalent to the category of coherent sheaves on a `noncommutative projective scheme', where a noncommutative projective scheme in the sense of Artin and Zhang is a quotient $\mathsf{QGr}A$ of the category of graded modules over a noncommutative $\Z$-graded algebra $A$ by the Serre subcategory of (graded) torsion modules. In \cite{artin_noncommutative_1994} the algebra $A$ is assumed to be noetherian, with this hypothesis weakened in \cite{polishchuk_noncommutative_2005} to include coherent algebras. Also in \cite{polishchuk_noncommutative_2005}, $\Z$-graded algebras are replaced by $\Z$-indexed algebras; a larger class of algebras (strictly containing $\Z$-graded algebras) which arise naturally in noncommutative projective geometry.

There are many examples where a $\Z$-grading/indexing is non-natural; perhaps the most notable of which being multi-projective spaces. These objects appear across both the commutative and noncommutative settings \cites{simis_diagonal_1998,chan_twisted_2000,presotto_birational_2018, presotto__2017}, and are graded/indexed by $\Z^r$ for some $r\geq 1$. Commonly we have $r=2$, as tensor products of $\Z$-graded/indexed algebras and blow-up algebras are both naturally $\Z^2$-graded/indexed. There are also examples beyond free abelian groups, such as the homogeneous coordinate rings of weighted projective lines \cite{geigle_class_1987}, which are graded by a finitely generated abelian group of rank 1, or the $\Z \times \{0,1,2\}$-indexed algebras arising from exceptional collections on noncommutative cubic surfaces in \cite{abdelgadir_cubic_2024}. 

When $I \neq \Z$, it is common to restrict an $I$-indexed algebra to some $\QGr$-equivalent $\Z$-indexed subalgebra (such as a Veronese or diagonal subalgebra) and work in the $\Z$-indexed setting. This is often undesirable; introducing arbitrary choice (for instance, the choice of an order preserving inclusion $\Z \hookrightarrow I$) or losing information from the original $I$-indexed algebra. The goal of this paper is to circumvent these issues by developing the theory of coherent $I$-indexed algebras and noncommutative projective schemes which unifies all of the aforementioned examples. Our main result (Theorem 5.10) is a version of Polishchuk's theorem when the index set $I$ is a locally finite directed poset, characterising the noncommutative projective spaces associated to coherent $I$-indexed algebras. This allows one to work with very general indexed algebras (including all the examples above) without being forced to choose some bijection between a natural index set and $\Z$, or restrict a larger indexed algebra to a $\Z$-indexed subalgebra. The trade-off for the level of generality of our index sets is that we must impose an extra condition, \textit{tails-finite generation} (c.f. Proposition 3.9), on our indexed algebras. This condition allows us to recover a number of results from the $\Z$-indexed setting which rely on the facts that $\Z$ is both totally ordered and a group, while requiring our index sets to be neither. 

Section 2 covers the basic theory of $I$-indexed algebras and their graded modules. Section 3 introduces the finiteness conditions we will need to impose on our indexed algebras and modules, as well as the construction of the noncommutative projective scheme $\qgr A$ associated to an $I$-indexed algebra $A$. In Section 3.2 the notion of tails-finite generation is introduced, ensuring that the subcategory of torsion modules (Definition 3.16) is localising (Theorem 3.18). Section 4 contains results concerning the cohomology of $I$-indexed algebras and noncommutative projective schemes necessary to prove the main result. In particular, we need an $I$-indexed version of the $\chi$ condition from \cite{artin_noncommutative_1994}. In Section 5, ample $I$-indexed sequences in abelian categories are introduced (Definition 5.1), which are analogous to ample $\Z$-indexed sequences as in \cite{polishchuk_noncommutative_2005}. We then show that ample $I$-indexed sequences characterise these noncommutative projective schemes up to equivalence, obtaining an $I$-indexed version of \cite[Theorem 2.4]{polishchuk_noncommutative_2005}. In Section 6 we show that the examples mentioned earlier are tails-finitely generated $I$-indexed algebras, and hence that our notion of noncommutative projective schemes encompasses these examples. Finally, Section 7 contains an application of the main theorem to the representation theory of bimodule species. Using Theorem 5.10 we show that the preprojective algebra $\Pi$ of a (possibly non-algebraic) bimodule species of infinite representation type is coherent, and that the corresponding noncommutative projective scheme $\qgr \Pi$ is derived equivalent to the category of finitely presented modules over the path algebra of the species (Theorem 7.3, Corollary 7.5). We view this as a version of Beilinson's derived equivalence \cite{beilinson_coherent_1978}, extending known results for hereditary finite dimensional algebras \cite{lenzing_curve_1986,baer_tame_1987,minamoto_ampleness_2012,krause_serre_2024} and for species over the Kronecker quivers \cite{chan_species_2016, chan_representation_2019}.

\section{$I$-Indexed algebras and modules}
We begin by discussing indexed algebras, also referred to as ringoids or algebroids elsewhere in the literature.

\begin{defn}{\cite[Section 2]{van_den_bergh_noncommutative_2011}}
Let $I$ be a set. An $I$\textbf{-indexed algebra} is a preadditive category $\A$ whose objects $\{ \A_i \}_{i \in I}$ are indexed by $I$.
\end{defn}

From an $I$-indexed algebra we may construct the ring $A = \bigoplus_{i,j \in I} A_{ij}$, where $A_{ij} = \Hom_{\A}(\A_j,\A_i)$. The addition on $A$ comes from the abelian groups $A_{ij}$, and multiplication is given by composition $A_{id} \times A_{dj} \to A_{ij}$. Note that $A$ is non-unital if $I$ is not finite.

\begin{ex}
    Any ring $R$ determines a $\{\ast\}$-indexed algebra $\mathbf{B}R$ with\break $\End_{\mathbf{B}R}(\ast) = R$. This is analogous to the construction of the delooping groupoid of a group.
\end{ex}

\begin{ex}{\cite{sierra_g-algebras_2011}}
    Let $G$ be a group and $S$ be a $G$-graded ring. Then we may form a $G$-indexed algebra $\mathbf{B}_G S$ as the full subcategory of $\Gr_G S$ consisting of the shifted modules $\{ S[g^{-1}] ~|~ g \in G \}$. As $\Gr_G S$ is always abelian, $\mathbf{B}_G S$ is preadditive. The associated ring is $\mathbf{B}_G S = \bigoplus_{g,h \in G} (\mathbf{B}_G S)_{g,h}$, where $(\mathbf{B}_G S)_{g,h} = S_{g^{-1}h}$ and multiplication is given by
    \[(\mathbf{B}_G S)_{g,h} \times (\mathbf{B}_G S)_{h,l} = S_{g^{-1}h} S_{h^{-1}l} \subset S_{g^{-1}l} = (\mathbf{B}_G S)_{g,l}.\]
    We note that Example 2.2 is a specific case of this construction when $G$ is the trivial group.
\end{ex}

We will frequently equate the ring $A$ with the $I$-indexed algebra $\A$. 

\begin{defn}{\cite[Section 2]{van_den_bergh_noncommutative_2011}}
    Let $\A$ be an $I$-indexed algebra. A \textbf{graded right} $\A$-\textbf{module} $M$ is a covariant additive functor $\A \to \Ab$.
\end{defn}

Equivalently, a graded right $A$-module is an abelian group $M = \bigoplus_{i \in I}M_i$ with abelian group homomorphisms $M_i \times A_{ij} \to M_j$. We will denote the category of graded right $A$-modules by $\Gr A$. Graded left modules are defined similarly via additive functors $\A^{op} \to \Ab$. We will usually refer to graded (left or right) $A$-modules simply as $A$-modules. Being defined in terms of functors to $\Ab$, we see immediately that $\Gr A$ is a Grothendieck category. 

\begin{ex}
    For a group $G$ and $S = \bigoplus_{g \in G} S_g$ a $G$-graded ring, $\Gr_G S \simeq \Gr (\mathbf{B}_GS)$. In fact, two $G$-graded rings $S$ and $S'$ are graded Morita equivalent if and only if $\mathbf{B}_GS \cong \mathbf{B}_GS'$ as $G$-indexed algebras \cite[Theorem 3.6]{sierra_g-algebras_2011}.
\end{ex}

While $A$ is generally non-unital, it always contains `local units' $e_i \in A_{ii}$ given by the identity morphisms $1_{\A_i}$. The local units determine right modules $e_iA$ which form a set of projective generators for $\Gr A$ \cite[Section 3]{van_den_bergh_noncommutative_2011}.

Given an $I$-indexed algebra $A$ and a $J$-indexed algebra $B$, we may form their \textbf{(external) tensor product}, defined by
\[A \otimes B := \bigoplus_{(i,j),(i',j') \in I \times J} A_{ij} \otimes_{\Z} B_{i'j'}.\]
Here external refers to the fact that these objects are no longer $I$-indexed nor $J$-indexed, but $I \times J$-indexed.

\begin{defn}{\cite[Section 2]{van_den_bergh_noncommutative_2011}}
    Let $\A$ be an $I$-indexed algebra and $\BB$ a $J$-indexed algebra. The category of $\A-\BB$-\textbf{bimodules} $\mathsf{Bimod}(\A-\BB)$ is the category of additive functors $[\A^{op} \otimes \BB, \Ab] = \Gr (\A^{op} \otimes \BB)$.
\end{defn}

Equivalently, an $A-B$-bimodule is a bigraded abelian group $\bigoplus_{(i,j) \in I \times J} T_{ij}$ with commuting homomorphisms $\mu_{ijj'}: T_{ij} \times B_{j,j'} \to T_{ij'}$ satisfying the right $B$-module axioms and $\psi_{i'ij}: A_{i',i} \times T_{ij} \to T_{i'j}$ satisfying the left $A$-module axioms. 

The units $e_i$ and $f_j$ of the $A_{ii}$ and the $B_{jj}$ give us restriction maps \[\mathsf{Bimod}(A-B) \to \Gr{A^{op}}\] and  \[\mathsf{Bimod}(A-B) \to \Gr{B},\]
given by $T \mapsto Tf_j$ and $T \mapsto e_iT$, respectively, for $T \in \mathsf{Bimod}(A-B)$.

An $A-B$-bimodule $T$ determines an adjoint pair of functors
$- \underline{\otimes}_A T: \Gr A \to \Gr B$ and $\underline{\Hom}_{B}(T_B,-): \Gr B \to \Gr A$, defined below in a similar fashion to the $\Z$-indexed case \cites{nyman_abstract_2019,mori_local_2021}. 

\begin{defn}
    Let $A$ and $B$ be $I$-indexed algebras, let $M \in \Gr A$, $N \in \Gr B$, and let $T$ be an $A-B$-bimodule. The \textbf{internal tensor product }of $M$ and $T$ is the right $B$-module
    \[M \underline{\otimes}_A T := \coker\left( \bigoplus_{l,m \in I} Me_l \otimes_{A_{ll}} A_{lm} \otimes_{A_{mm}} e_mT \xrightarrow{\mu \otimes 1 - 1 \otimes \mu} \bigoplus_{n \in I} Me_n \otimes_{A_{nn}} e_n T \right),\]
    where $\mu$ denotes scalar $A$-multiplication on a (left or right) $A$-module.
    The \textbf{internal hom} of $T$ and $N$ is the right $A$-module
    \[\underline{\Hom}_B(T,N) := \bigoplus_{i \in I} \Hom(e_iT,N),\]
    with the right $A$-module structure coming from the left $A$-module structure on $T$.
\end{defn}

Viewing $A$ as an $A$-bimodule, we have that $\underline{\Hom}_A(A,M) \cong M \cong M \underline{\otimes}_A A$ for any $M \in \Gr A$. 

\section{Coherent $I$-indexed algebras and $\qgr$}
The aim of this section is to construct the noncommutative projective scheme $\qgr A$ associated to a coherent $I$-indexed algebra $A$. For the remainder of the paper we will suppose that $I$ is a \textbf{locally finite poset}: a partially ordered set $(I,\leq)$ such that for any $i,j \in I$, the interval $[i,j] := \{ d \in I ~|~ i \leq d \leq j \}$ is finite. We will also use the notation $(i,j) := \{ d \in I ~|~ i < d < j \}$ for open intervals. If $i,j \in I$ are incomparable, so that neither $i<j$ nor $j<i$, we write $i\vert\vert j$.

\begin{ex}
    \begin{enumerate}
        \item Let $I = (\Z^r,\leq_{prod})$ be $\Z^r$ endowed with the product order. For any $(m)=(m_1,\ldots,m_r),~(n)=(n_1,\ldots,n_r) \in \Z^r$, $[(m),(n)]$ is finite as it is simply the product of the finite intervals $[m_i,n_i]$ in $\Z$ for $1 \leq i \leq r$.
        \item Let $G$ be a finitely generated abelian group of rank 1, in which case there are positive integers $w_1,\ldots,w_n$ such that $G \cong \Z \oplus \Z/w_1\Z \oplus \ldots \oplus \Z/w_n\Z$. We may choose generators $x_1,\ldots,x_n$ for $G$ such that the relations are given by $w_1x_1=\ldots=w_nx_n=:c$. Then $G$ is naturally partially ordered by setting $mc + \sum_{l=1}^n m_lx_l \leq_G m'c + \sum_{l=1}^n m'_{l}x_{l}$ if $m<m'$ or $m=m'$ and $m_l\leq m'_l$ for $l=1,\ldots,n$. It is clear that $(G,\leq_G)$ is a locally finite poset as $\Z$ is and the $\Z/w_l\Z$ are each finite.
        \item Let $Q=(Q_0,Q_1)$ be a quiver which is finite, connected and acyclic. From $Q$ one constructs a \textit{translation quiver} $\Z Q$, which plays a fundamental role in the representation theory of $Q$, and of finite dimensional algebras over a field more generally \cites{auslander_artin_1995, happel_triangulated_1988}. This quiver has vertex set $(\Z Q)_0 = \Z \times Q_0$ and edge set $(\Z Q)_1$ consisting of edges $e_m:(m,v) \to (m,w)$ and $e_m^{\ast}:(m,w) \to (m+1,v)$ for each edge $e:v \to w \in Q_1$ and every $m \in \Z$. Now $\Z Q$ is locally finite and acyclic as $Q$ is finite and acyclic, so $\Z Q$ determines a locally finite partial order $(\Z \times Q_0,\prec)$ where $(m,v) \prec (n,w)$ if there is a directed path from $(m,v)$ to $(n,w)$ in $\Z Q$.
    \end{enumerate}
\end{ex}

Given an $I$-indexed algebra $A$, we let $A_{\Delta} := \bigoplus_{i \in I} A_{ii}$. Then $\diag{A}$ is naturally a subalgebra of $A$ which we call the \textbf{base} of $A$. The base $\diag{A}$ is analogous to the degree 0 component of a graded ring, so we say $A$ is \textbf{connected} if each $A_{ii}$ is a division ring. We say an $I$-indexed algebra is \textbf{positively indexed} (positively $I$-indexed) if $A_{ij} \neq 0$ only when $i \leq j$. 

If $A$ is connected we will usually write $A_{ii} = K_i$. When $A$ is connected and positively indexed, $\m = \bigoplus_{i<j}A_{ij}$ is the unique maximal ideal of $A$. The quotients $S_i := e_iA/ e_i\m$ for $i \in I$ are then simple right $A$-modules, and moreover any simple right $A$-module is isomorphic to some $S_i$. We also note that when $A$ is connected, multiplication on $A$ gives us, for all $i,d,j \in I$, $K_i-K_j$-bimodule homomorphisms $\mu_{idj}:A_{id} \otimes_{K_d} A_{dj} \to A_{ij}$ defined by $\mu_{idj}(\sum_l x_l \otimes y_l) = \sum_l x_ly_l$.

\subsection{Finiteness conditions for $I$-indexed algebras}
For the remainder of the paper we will suppose that any $I$-indexed algebra $A$ is connected and positively indexed, with base $\K = \diag{A} = \bigoplus_{i \in I} K_i$. We will work with right modules unless otherwise specified. Recall that $\{P_i := e_iA\}_{i \in I}$ is a set of projective generators for $\Gr A$. We say a right module $F$ is \textbf{free} if it is a direct sum of the $P_i$. If this sum is finite, we say $F$ is \textbf{free and finitely generated}. 

\begin{defn}
    An $A$-module $M$ is \begin{itemize}
        \item \textbf{finitely generated} if there is a surjection $F \twoheadrightarrow M$ from some free and finitely generated module $F$, 
        \item \textbf{finitely presented} if there is an exact sequence $F' \to F \to M \to 0$, with $F,F'$ free and finitely generated, 
        \item \textbf{locally finite} if $\dim_{K_j}M_j$ is finite for all $j \in I$, 
        \item \textbf{finite dimensional} if $M$ is locally finite and $M_i \neq 0$ for finitely many $i$.
    \end{itemize}
\end{defn}
We will also call an $I$-indexed algebra $A$ (right) \textbf{locally finite} when $A_{ij}$ is finite dimensional over $K_j$ for all $i,j \in I$. We will refer to a right locally finite $I$-indexed simply as locally finite, as we work only with right modules, although left and two-sided local finiteness may be defined similarly.

\begin{lemma}
    $A$ is locally finite if and only if $P_i$ is locally finite for all $i \in I$.
\end{lemma}
\begin{proof}
    This is immediate as $A_{ij} = e_iAe_j =  (P_i)_j$.
\end{proof}

\begin{defn}
Let $A$ be an $I$-indexed algebra. An $A$-module $M$ is \textbf{coherent} if $M$ is finitely generated and the kernel of any homomorphism $F \to M$, where $F$ is free and finitely generated, is also finitely generated.
\end{defn}

We now let $\gr A$ denote the full subcategory of finitely presented $A$-modules, and let $\mathsf{coh} A$ be the full subcategory of coherent $A$-modules. The latter is an abelian subcategory of $\mathsf{Gr} A$ which is closed under extensions. Recall (\cite{herzog_ziegler_1997}) that an object $X$ in a Grothendieck category is \begin{itemize}
    \item \textbf{finitely generated} if for any collection of subobjects $\{X_j \subset X \}_{j \in J}$ such that $\sum_{j \in J} X_j = X$, there is a finite subset $J' \subset J$ with $X = \sum_{j \in J'}X_j$. 
    \item \textbf{coherent} if $X$ is finitely generated and for any morphism $f:Y \to X$ with $Y$ finitely generated, $\ker f$ is finitely generated.
\end{itemize}
Of course, these definitions agree with the above definitions for finitely generated and coherent modules when the Grothendieck category is $\Gr A$. We also recall that a Grothendieck category is \textbf{locally coherent} if it possesses a family of coherent generators, or equivalently every object is a directed colimit of coherent objects. 

\begin{propn}
    Let $A$ be an $I$-indexed algebra. The following are equivalent:
    \begin{enumerate}
        \item $\gr A = \coh A$, and as such is a full abelian subcategory of $\Gr{A}$, 
        \item $\Gr A$ is a locally coherent category, 
        \item The free right modules $P_i$ are coherent for all $i \in I$.
    \end{enumerate}
\end{propn}
\begin{proof}
    $(1) \Leftrightarrow (2)$ is well known (see \cite[Section 2]{roos_locally_1969} or \cite[Theorem 1.6]{herzog_ziegler_1997}).
    
    $(1) \Rightarrow (3)$ follows as the $P_i$ are finitely presented. Finally, every finitely presented module is a cokernel of a morphism between finite sums of the $P_i$, so if all the $P_i$ are coherent the cokernel of this morphism must be too. Hence $(3) \Rightarrow (1)$.
\end{proof}

If any of the above equivalent conditions are satisfied, we say $A$ is a (right) \textbf{coherent} $I$-indexed algebra. If, in addition, all the simple modules $S_i$ are coherent, we say $A$ is \textbf{strongly coherent}. We denote the full subcategory of finite dimensional modules by $\mathsf{fd}A$. As any simple $A$-module is isomorphic to some $S_i$, we see that $\mathsf{fd}A \subset \gr A$ if $A$ is strongly coherent.

\begin{remark}
    If we let $\ssf{fg}A$ denote the full subcategory of finitely generated modules, and let $\ssf{ACC}(A)$ denote the full subcategory of noetherian modules (modules satisfying the ascending chain condition on submodules) then, similarly to Proposition 3.5, we have that $\ssf{fg}A = \ssf{ACC}(A)$ if and only if $\Gr A$ is locally noetherian ($\Gr A$ has a family of generators satisfying the ascending chain condition on subobjects) if and only if $P_i$ is noetherian for all $i \in I$. An $I$-indexed algebra satisfying any of these equivalent conditions is called (right) \textbf{noetherian}. Any quotient of a noetherian module is again noetherian, however, so there is no need to define a notion of `strong noetherian' $I$-indexed algebras.
\end{remark}

\subsection{Tails of graded modules}
A submodule of $M \in \Gr A$ of the form $M_{\geq d} := \bigoplus_{j \geq d} M_j$, for some $d \in I$, is called a \textbf{(weak) tail} of $M$. Similarly, we call $M_{>d}:= \bigoplus_{j > d} M_j$ a \textbf{(strict) tail} of $M$. For any $A$-module $M$ and $d \in I$, the strict tail $M_{>d}$ embeds in the weak tail $M_{\geq d}$, so we have exact sequences
\begin{equation}
    0 \to M_{>d} \to M_{\geq d} \to M_d \to 0.
\end{equation}
In this way, each $M_d$ has the structure of a right $A$-module. We also note that we have exact sequences
\begin{equation}
    0 \to M_{>d} \to M \to M/M_{>d} \to 0
\end{equation}
for all $M \in \Gr A$ and $d \in I$.

When $I=\Z$, $M_{>n} = M_{\geq n+1}$ for all $n \in \Z$, so every strict tail is a weak tail (with the index shifted). This may no longer be true if $I\neq \Z$. In most cases we will work with strict tails, referring to them simply as tails.

We now introduce two conditions on an $I$-indexed algebra which are equivalent when $I=\Z$ \cite[Proposition 3.2]{de_deken_abelian_2011}. 
\begin{enumerate}
    \item[($\ast$)] For all $i \in I$, $e_i \m = P_{i,>i}$ is finitely generated. 
    \item[($\ast\ast$)] For all $i,d \in I$, $P_{i,>d}$ is finitely generated.
\end{enumerate}
If $I$ is not totally ordered, this equivalence need not hold.

\begin{ex}
    Let $k$ be a field and consider the free algebra $k\langle x,y \rangle$ as a $\Z^2$-graded algebra, with $x$ in degree $(1,0)$ and $y$ in degree $(0,1)$. If $\Z^2$ is given the product order, the associated $\Z^2$-indexed algebra $B=\mathbf{B}_{\Z^2}k\langle x,y \rangle$ satisfies ($\ast$), as
    \[B_{(m,n),>(m,n)} = xB_{(m+1,n)} + yB_{(m,n+1)}\]
    for all $(m,n) \in \Z^2$. The tail $e_{(0,0)}B_{> (1,1)}$ is not a finitely generated module, however, as the collection of elements $\{x^2y,x^3y,x^4y, \ldots \} \subset e_{(0,0)}B_{> (1,1)}$ is pairwise (right) linearly independent, and so ($\ast\ast$) does not hold.
\end{ex}

While ($\ast$) may not imply ($\ast\ast$), it is still sufficient to ensure $A$ is locally finite, similarly to \cite[Lemma 3.1]{de_deken_abelian_2011}.

\begin{lemma}
    If $A$ satisfies ($\ast$), then $A$ is locally finite.
\end{lemma}
\begin{proof}
    We show $A_{ij}$ is finite over $K_j$ for all $i,j \in I$ by inducting on $v_{ij} := \vert [i,j] \vert$. The case $v_{ij} =1$ is trivial, as we must have $i=j$. Now suppose $v_{ij}>1$, and that $A_{i'j'}$ has finite dimension over $K_{j'}$ for all $i',j' \in I$ with $v_{i'j'} < v_{ij}$. Multiplication gives us a map $\mu_{(i,j)} := \bigoplus_{d \in (i,j)} \mu_{idj}:\bigoplus_{d \in (i,j)} A_{id} \otimes_{K_d} A_{dj} \to A_{ij}$, so we have an exact sequence of $K_i-K_j$-bimodules
    \begin{equation}
        0 \to \mathrm{im}~\mu_{(i,j)} \to A_{ij} \to \coker \mu_{(i,j)} \to 0.
    \end{equation}
    The inductive hypothesis tells us that $A_{id}$ and $A_{dj}$ are finite over $K_d$ and $K_j$ respectively, as $v_{id} ,v_{dj} <v_{ij}$. Then both $\bigoplus_{d \in (i,j)} A_{id} \otimes_{K_d} A_{dj}$ and the quotient $\mathrm{im}~\mu_{(i,j)}$ are finite over $K_j$. Now take any $a \in A_{ij}$ with image $\overline{a}$ in $\coker \mu_{(i,j)}$. Using ($\ast$), for $i \in I$ we let $X_{i,>i} = \{ x_l \in A_{ij_l}\}_{l=1}^{r}$ be a finite generating set for $P_{i,>i}$. Then we have $a=\sum_{l=1}^{r} x_lb_l$ for some $b_l \in A_{j_lj}$. If $j_l \neq j$, $x_lb_l \in \mathrm{im}~\mu_{(i,j)}$, so   
    \[\overline{a} = \sum_{j_l =j} x_lb_l + \mathrm{im}~\mu_{(i,j)}.\]
    Then, as $X_{i,>i}$ is finite, $\{ x + \mathrm{im}~\mu_{(i,j)} \mid x \in X_{i,>i} \cap A_{ij} \}$ is a finite $K_j$-basis for $\coker\mu_{(i,j)}$. It now follows from (3.3) that $A_{ij}$ is finite over $K_j$.
\end{proof}

In \cite{de_deken_abelian_2011}, a $\Z$-indexed algebra satisfying $(\ast)$ (or equivalently $(\ast \ast)$) is called \textit{finitely generated}.

\begin{propn}
    For an $I$-indexed algebra $A$, the following are equivalent. 
    \begin{enumerate}
        \item $A$ satisfies $(\ast\ast)$. 
        \item $P_i / P_{i,>d}$ is finitely presented for all $i,d \in I$. 
        \item $M_{>d}$ is finitely generated for all $M \in \gr A$ and all $d \in I$. 
        \item $M/M_{>d}$ is finitely presented for all $M \in \gr A$ and all $d \in I$.
    \end{enumerate}
    If $A$ satisfies any of the equivalent conditions above, we say $A$ is (right) \textbf{tails-finitely generated}.
\end{propn}
\begin{proof}
    $(1) \Leftrightarrow (2)$ follows from the exact sequence (3.2) as $P_i$ is finitely presented. $(3) \Rightarrow (1)$ is trivial, and $(1) \Rightarrow (3)$ follows by taking the tails of a finite presentation of $M$. Finally, $(3) \Leftrightarrow (4)$ also follows from (3.2) as $M \in \gr A$.
\end{proof}

If $A$ is tails-finitely generated then the strict and weak tails of any $M \in \gr A$ must be finitely generated.

\begin{propn}
    If $A$ is tails-finitely generated and $M \in \gr A$, then the $A$-modules $M_{\geq d}$ and $M_d$ are finitely generated for all $d \in I$.
\end{propn}
\begin{proof}
    It suffices to show for $M=P_i$ by considering a finite presentation of $M$. Lemma 3.8 says $(P_i)_d = A_{id}$ is finitely generated, and so tails-finite generation and (3.1) ensure $P_{i,\geq d}$ is also finitely generated.
\end{proof} 

\begin{ex}
    Once again we let $k$ be a field, but now we consider the commutative polynomial ring $R= \bigoplus_{(m,n) \in \Z^2} R_{(m,n)} = k[x,y]$, graded by $\deg(x)=(1,0)$ and $\deg(y)=(0,1)$. Then $D=\mathbf{B}_{\Z^2}R$ satisfies ($\ast$) similarly to Example 3.7, however $D$ is also tails-finitely generated. For $(m,n),(r,s) \in \Z^2$, commutativity shows 
    \[e_{(m,n)}D_{> (r,s)} \cong D_{(m,n)(r+1,s)} (e_{(r+1,s)}D) + D_{(m,n)(r,s+1)} (e_{(r,s+1)}D),\]
    and both $D_{(m,n)(r+1,s)} \cong R_{(r-m+1,s-n)}$ and $D_{(m,n)(r,s+1)} \cong R_{(r-m,s-n+1)}$ are finite dimensional over $K_{(r+1,s)} \cong K_{(r,s+1)} \cong k$. It then follows that $e_{(m,n)}D_{>(r,s)}$ is finitely generated, so $D$ is tails-finitely generated.
\end{ex}

We present some more examples of tails-finitely generated indexed algebras in Section 6. For now, we look at some situations in which an $I$-indexed algebra is tails-finitely generated. Given two elements $i,d \in I$, we define a set $U_{min}(i,d) = \{ j \in [i,-) \cap [d,-) \mid [i,j) \cap [d,j) = \varnothing \}$ containing all of the minimal upper bounds for $i$ and $d$. We also recall the multiplication maps $\mu_{idj}: A_{id} \otimes_{K_d} A_{dj} \to A_{ij}$ for $i,j,d \in I$ defined earlier.

\begin{propn}
    Let $A$ be an $I$-indexed algebra.
    \begin{enumerate}
        \item $A$ is tails-finitely generated if $A$ is noetherian. 
        \item $A$ is tails-finitely generated if $A$ satisfies ($\ast$) and either \begin{enumerate}
            \item $\vert U_{min}(i,d) \vert < \infty$ and $\mu_{iju}:A_{ij} \otimes_{K_j} A_{ju} \to A_{iu}$ is surjective for all $i,d,j,u \in I$ with $i \leq j \leq u$, or
            \item $\vert \{j \in I \mid j \geq i,~ j\vert\vert d\} \vert < \infty$ for all $i,d \in I$.
        \end{enumerate}
    \end{enumerate}
\end{propn} 
\begin{proof}
    (1): If $A$ is noetherian all the $P_{i}$ are noetherian modules (c.f. Remark 3.6), so the submodules $P_{i,>d}$ are all finitely generated.

    (2a): We consider three cases. If $d \leq i$ then, as $A$ is positively indexed, $P_{i,>d}$ is either $P_i$ or $P_{i,>i}$. The free module $P_i$ is generated by $e_i \in (P_i)_i$, and $P_{i,>i}$ is finitely generated as $A$ satisfies $(\ast)$.

    Suppose $i<d$. Using $(\ast)$ and Lemma 3.8 we take a finite $K_d$-basis $X$ for $A_{id}$ and a finite generating set $Y$ for $P_{d,>d}$. If $(P_{i,>d})_u \neq 0$ for $u \in I$, then $(P_{i,>d})_u = A_{iu}$. As the multiplication maps are all surjective, for any $0 \neq m \in (P_{i,>d})_u$ there is some $\sum_l a_l \otimes b_l \in A_{id} \otimes_{K_d} A_{du}$ such that $m = \mu_{idu}(\sum_l a_l \otimes b_l) = \sum_l a_lb_l$. Using the basis $X$ we may write $a_l = \sum_{x \in X} xc_{l,x}$ for some $c_{l,x} \in K_d$, so that $m = \sum_l \sum_{x \in X} xc_{l,x}b_l$. Now each $c_{l,x}b_l \in A_{du}= (P_{d,>d})_u$, so using the generating set $Y$ we may write $c_{l,x}b_l = \sum_{y \in Y} yc'_{l,x,y}$ for some $c'_{l,x,y} \in Ae_u$. If we let $c''_{x,y} = \sum_l c'_{l,x,y}$, we then have
    \[ m = \sum_l \sum_{x \in X} xc_{l,x}b_l = \sum_l \sum_{x \in X} \sum_{y \in Y} xyc'_{l,x,y} = \sum_{\substack{x \in X \\ y \in Y}}xy(c''_{x,y}) \in \sum_{\substack{x \in X \\ y \in Y}} xy A. \]
    Thus $\{ xy \mid x \in X,~y \in Y \}$ is a finite generating set for $P_{i,>d}$.

    Finally, we suppose $d \vert \vert i$. For each $j \in U_{min}(i,d)$ we use Lemma 3.8 to take a finite $K_j$-basis $X_j$ for $A_{ij}$. Now take some $0 \neq m \in (P_{i,>d})_u$. As $A$ is positively indexed, $(P_{i,>d})_u \neq 0$ only if $u \in [i,-) \cap (d,-)$. Thus there is some $j \in U_{min}(i,d)$ such that $j \leq u$. The multiplication map $\mu_{iju}:A_{ij} \otimes_{K_j} A_{ju} \to A_{iu}=(P_{i,>d})_u$ is surjective by assumption, so there is some $\sum_l a_l \otimes b_l \in A_{ij} \otimes_{K_j} A_{ju}$ such that $m = \mu_{iju}(\sum_l a_l \otimes b_l) = \sum_l a_lb_l$. Using the basis $X_j$ we may write $a_l = \sum_{x \in X_j} xc_{l,x}$ for some $c_{l,x} \in K_j$, so that we have
    \[ m = \sum_l a_lb_l = \sum_l \sum_{x \in X_j} xc_{l,x}b_l = \sum_{x \in X_j} x(\sum_l c_{l,x}b_l) \in \sum_{x \in X_j} x A. \]
    It follows that $X = \bigcup_{j \in U_{min}(i,d)} X_j$ is a generating set for $P_{i,>d}$. Moreover, $X$ is finite as each $X_j$ is finite and $U_{min}(i,d)$ is finite by assumption, so $P_{i,>d}$ is finitely generated.

    (2b): Suppose $i,d \in I$. We note again that if $d<i$ then $P_{i,>d} = P_i$ is always finitely generated, so we suppose $d \not< i$. Let $J_{id} = \{ j \in I \mid j \geq i, j \not> d \} = [i,d] \cup \{ j \in I \mid j \geq i, j \vert \vert d \}$. As $I$ is locally finite and $\{ j \in I \mid j \geq i, j \vert \vert d \}$ is finite by assumption, we have that $J_{id}$ is finite. We induct on $\vert J_{id}\vert$. Note that $i \in J_{id}$, so when $\vert J_{id} \vert = 1$ we must have $i=d$. In this case $P_{i,>i}$ is finitely generated by $(\ast)$. Now suppose $\vert J_{id} \vert = n$ and that $P_{i',>d'}$ is finitely generated for all $i',d' \in I$ with $\vert J_{i'd'} \vert < n$. Using $(\ast)$ we take a finite generating set $\{ x_l \in A_{i,j_l}\}_{l=1}^r$ for $P_{i,>i}$. As $i < j_l$ for all $l$, we have $\vert J_{j_ld} \vert < \vert J_{id} \vert$ for all $l$. Hence each $P_{j_l,>d}$ is finitely generated by the inductive hypothesis. So if $Y_l$ is a finite generating set for $P_{j_l,>d}$ for $1 \leq l \leq r$, then $\{ x_ly \mid 1 \leq l \leq r,~y \in Y_l\}$ is a finite generating set for $P_{i,>d}$.
    \end{proof}

Tails-finite generation implies the simple modules $S_i = P_i / P_{i,>i}$ are finitely presented, so the following is immediate.

\begin{lemma}
    If $A$ is tails-finitely generated, then $A$ is strongly coherent if and only if $A$ is coherent.
\end{lemma}

We will refer to an $I$-indexed algebra that is tails-finitely generated and coherent as \textbf{tails-coherent}, as it follows from Propositions 3.9 and 3.10 that $A$ is tails-coherent if and only if $A$ is coherent and every (strict or weak) tail of a coherent module is coherent. 

\subsection{$\qgr$ for tails-coherent $I$-indexed algebras}
For the remainder of the paper all $I$-indexed algebras will be tails-finitely generated and, in addition to being a locally finite poset, we insist $I$ is \textbf{directed} so that for any $i,j \in I$, there is some $u_{ij} \in I$ with $i,j \leq u_{ij}$. 

\begin{ex}
    \begin{enumerate}
        \item Let $(\Z^r,\leq_{prod})$ be $\Z^r$ endowed with the product order. Any $(m),(n) \in \Z^r$ have a least upper bound $(\max\{m_1,n_1\},\ldots,\max\{m_r,n_r\}) \in \Z^r$, so $(\Z^r,\leq_{prod})$ is directed. 
        \item Let $(G,\leq_G)$ be the partial order on a finitely generated abelian group $G$ of rank 1 described in Example 3.1(2). Two elements $mc + \sum_{l=1}^n m_lx_l$ and $m'c + \sum_{l=1}^n m'_lx_l$ in $G$ are upper-bounded by $(\max\{m,m'\}+1)c$, so $(G,\leq_G)$ is directed.
        \item Let $I =(\Z \times Q_0,\prec)$ be the partial order defined by the translation quiver $\Z Q$ as in Example 3.1(3). Given an undirected path (path in the undirected graph underlying $Q$) $p=e_1e_2\ldots e_m$ in $Q$, we let $r(p)$ denote the number of edges in $p$ whose orientation differs from the orientation of $Q$. In this way, $p$ is a directed path in $Q$ if and only if $r(p)=0$. Then we may define the partial order on $\Z \times Q_0$ by stating $(m,v) \prec (n,w)$ if there is an undirected path $p:v \to w$ in $Q$ with $r(p) \leq n-m$. As $Q$ is finite and connected, there is an undirected path $p_{vw}:v \to w$ for all $v,w \in Q_0$ with $r(p) \leq \vert Q_0 \vert$. It follows that $(\Z \times Q_0,\prec)$ is directed.
    \end{enumerate}
\end{ex}

\begin{remark}
    It is not difficult to see that if $(I,\leq)$ is one of the posets in Example 3.14(2) or (3), any set of the form $\{j \in I \mid j \geq i,~ j\vert\vert d\}$ for $i,d \in I$ must be finite, and thus Proposition 3.12(2b) shows $(\ast)$ is equivalent to $(\ast \ast)$ for any algebra positively indexed by these partial orders.
\end{remark}

For an $I$-indexed algebra $A$ and an $A$-module $M$, we say $M$ is \textbf{upper-bounded} if $M_{> d}=0$ for some $d \in I$. 

\begin{defn}
    An $A$-module $M$ is \textbf{torsion} if it is a directed colimit of upper-bounded modules.
\end{defn}

Clearly $M$ is torsion if and only if the cyclic module $mA$ is upper-bounded for every $m \in M$. An element $m \in M$ with $mA$ upper-bounded is called a \textbf{torsion element}. If $M$ has no non-zero torsion elements, we say $M$ is \textbf{torsion-free}. We call the full subcategory of torsion modules $\mathsf{Tors} A$. Any $A$-module $M$ contains a largest torsion submodule $\tau(M) := \{ m \in M ~|~ m \text{ is torsion}\}$, and $M$ is torsion if and only if $\tau(M)=M$. The assignment $M \mapsto \tau(M)$ is functorial, and we call the functor $\tau:\Gr A \to \ssf{Tors}A$ the \textbf{torsion functor}.
 
For $\Z$-indexed algebras, the torsion functor has an explicit description 
\[\tau(-) \simeq \lim\limits_{n \to \infty} \HOM_A(A/A_{\geq n},-),\]
where  $A_{\geq n} = \bigoplus_{j-i \geq n} A_{ij}$ \cite[Lemma 6.7]{nyman_abstract_2019}. When $I\neq \Z$ we must define the torsion functor slightly differently. For $d \in I$, we define the two-sided ideal $A_{\ast,> d} = \bigoplus_{\substack{i,j \in I \\ j > d}} A_{ij}$ of $A$. We then have the exact sequence
\begin{equation}
    0 \to A_{\ast,> d} \to A \to A / A_{\ast,> d} \to 0
\end{equation}
in $\mathsf{Bimod}(A)$. 

\begin{theorem}
    Let $M \in \Gr A$. Then $\tau(M) \cong \varinjlim_{d \in I}\HOM_A(A/A_{\ast,> d},M)$.
\end{theorem}
\begin{proof}
    Applying $\varinjlim_d\HOM_A(-, M)$ to (3.4) shows $\varinjlim_{d \in I}\HOM_A(A/A_{\ast,> d},M)$ is a submodule of $\HOM_A(A,M) \cong M$. Each term $\HOM_A(A/A_{\ast,> d},M)$ is upper-bounded by $d$, so $\varinjlim_{d \in I}\HOM_A(A/A_{\ast,> d},M)$ is torsion and thus a submodule of $\tau(M)$. On the other hand, any $m \in \tau(M)_i$ is an element of $M_i\cong \Hom_A(P_i,M)$ such that $mA$ is upper-bounded. If the upper bound is $d \in I$, then $m$ is an element of $\Hom_A(P_i/P_{i,> d},M) = \HOM_A(A/A_{\ast,> d},M)_i$. So $m$ belongs to $\varinjlim_{d \in I}\HOM_A(A/A_{\ast,> d},M)$, and hence $\tau(M) \cong \varinjlim_{d \in I}\HOM_A(A/A_{\ast,> d},M)$. 
\end{proof}

\begin{theorem}
$\mathsf{Tors}A$ is a localising subcategory of finite type in $\Gr A$.
\end{theorem}
\begin{proof}
It follows from \cite[III.3 Corollaire 1]{gabriel_categories_1962} that $\mathsf{Tors}A$ is localising if and only if it is a Serre subcategory, as $\mathsf{Gr} A$ has enough injectives and Theorem 3.17 shows a torsion functor exists. One direction is easy, as any submodule or quotient of a torsion module is torsion.

Now suppose $0 \to M' \to M \to M'' \to 0$ is exact with both $M'$ and $M''$ torsion, and take any element $m \in M_i$. As $M''$ is torsion, there is some $d \in I$ such that $(mA)_{> d} \subseteq M'$. As $A$ is tails-finitely generated we may take a finite set of generators $\{ x_l \in A_{ij_l} \}_{l=1}^r$ for $P_{i,> d}$. Note that all of the $m x_{l}$ are torsion, as they are elements of $M'$. As there are only finitely many $x_l$ and $I$ is directed, there is some common upper bound $u \in I$ for $j_1,\ldots,j_r$, so that $(m x_{j_l}A)_{> u}= 0$ for $l=1,\ldots,r$. Then $mA$ is upper-bounded by $u$, so $M$ is torsion.

Finally, Proposition 3.9(2) and the description of $\tau$ in Theorem 3.17 show $\tau$ commutes with directed colimits, so $\ssf{Tors}A$ is of finite type.
\end{proof}

We let $\QGr A := \Gr A / \ssf{Tors}A$ denote the quotient category, which we view as the category of quasicoherent sheaves on the `noncommutative projectivisation' of $A$. We let $\pi: \mathsf{Gr} A \to \mathsf{QGr} A$ be the quotient functor, and as $\ssf{Tors} A$ is localising there is a section functor $\omega: \mathsf{QGr}A \to \mathsf{Gr}A$ which is right adjoint to $\pi$. The the next result follows from the previous theorem and the results of \cite[Section 2]{krause_coherent_1997}.

\begin{cor}
    Suppose $A$ is tails-coherent. Let $\iota$ be the inclusion of $\gr A$ in $\Gr A$, and let $\ssf{tors}A = \ssf{Tors}A \cap \gr A$.
    \begin{enumerate}
        \item $\ssf{tors}A$ is a Serre subcategory of $\gr A$, and $\ssf{Tors}A$ is the closure of $\ssf{tors}A$ under directed colimits. 
        \item The section functor $\omega$ commutes with directed colimits. 
        \item $\QGr A$ is locally coherent, and there is an equivalence $\Phi: \qgr A := \gr A / \ssf{tors} A \simeq \coh (\QGr A)$. 
        \item There is an exact, fully faithful inclusion $\underline{\iota}:\qgr A \to \QGr A$, given by composing $\Phi$ with the inclusion of $\coh(\QGr A)$ into $\QGr A$, such that $\pi_{\QGr} \iota = \underline{\iota}\pi_{\qgr}$. 
    \end{enumerate}
\end{cor}

Any $M \in \ssf{tors}A$ has a finite set of generators $\{m_1,\ldots,m_r\}$, and $m_lA$ is upper bounded for $1 \leq l \leq r$ as $M$ is torsion. As $I$ is directed we may choose a single upper bound for all of these cyclic modules, so $M$ is upper bounded. When $I=\Z$, we have $\ssf{tors}A = \ssf{fd}A$ \cite[Lemma 3.9]{mori_char_2025}. If $I$ is not totally ordered then a module which is finitely generated and upper bounded need not be finite dimensional (take the module $e_{(0,0)}B/e_{(0,0)}B_{>(1,1)}$ from Example 3.7, for example). If the index set $I$ satisfies the condition of Proposition 3.12(2b), however, then we recover this equality.

\begin{lemma}
    If $A$ is tails-coherent and the sets $\{j \in I \mid j \geq i,~ j\vert\vert d\}$ are finite for all $i,d \in I$, then $\ssf{tors}A=\ssf{fd}A$.
\end{lemma}
\begin{proof}
    As $A$ is tails-coherent, Lemma 3.8 ensures $A$ is locally finite, and hence it suffices to show that $M_j \neq 0$ for finitely many $j \in I$ for any $M \in \ssf{tors}A$. Suppose $M \in \ssf{tors}A$ is upper bounded by $d \in I$. If $\{ m_l \in M_{i_l} \}_{l=1}^r$ is a finite set of generators for $M$, then $M_j \neq 0$ only if $j \in J_M:=\bigcup_{l=1}^r \{ j \in I \mid j \geq i_l,~j \not> d \}$. As $\{ j \in I \mid j \geq i,~j \not> d \} = [i,d] \cup \{ j \in I \mid j \geq i,~ j \vert \vert d\}$ for all $i \in I$, we see that $J_M$ is a finite union of finite sets, and is thus finite. So $M_j$ is non-zero for only finitely many $j \in I$.
\end{proof}

In the future we will denote objects in the quotient category as $\MM = \pi M$, and will refer to such objects as (quasi-)coherent sheaves (coherent if $\MM \in \qgr A$). In particular, $\pi P_i$ is a coherent sheaf which we denote by $\cP_i$. We will use the notation $\QHom_A(-,-) := \Hom_{\QGr A}(-,-)$ for the morphisms in $\QGr A$. Corollary 3.19(4) identifies $\qgr A$ as a full subcategory of $\QGr A$, and so we will use the same notation for morphisms in $\qgr A$. For an $A$-bimodule $T$, we use $\QHOM_A(\mathcal{T},-)$ to denote $\bigoplus_{i \in I} \QHom_A(\pi (e_iT),-)$. 

\section{Cohomology of $I$-indexed algebras}
Throughout this section $A$ will be a tails-coherent $I$-indexed algebra. We collect some results about the cohomology of tails-coherent $I$-indexed algebras and the associated noncommutative projective schemes. Most of the results of this section are standard (see \cite[Section 4]{bondal_generators_2003} for the $\Z$-graded version, or \cite[Section 6]{nyman_abstract_2019} for the $\Z$-indexed version). The following is clear from the definition of coherence.

\begin{propn}
Every $M \in \gr A$ has a resolution of the form
\[ ... \to F^1 \to F^0 \to M \to 0,\]
where the $F^i$ are all free and finitely generated.
\end{propn}

The right derived functors of 
\[\underline{\Hom}_A(-,-): \mathsf{Bimod}(A)^{op} \times \Gr A \to \Gr A\]
are given by
\[\underline{\Ext}_A^l(-,-) := \bigoplus_{i \in I} \Ext_A^l(e_i(-),-):\mathsf{Bimod}(A)^{op} \times \Gr A \to \Gr A.\] 

\begin{lemma}
The right derived functors of the torsion functor $\tau: \mathsf{Gr}A \to \mathsf{Tors}A$ commute with directed colimits.
\end{lemma}
\begin{proof}
    Theorem 3.17 shows $\ssf{R}^l\tau(-) \cong \varinjlim_d\underline{\Ext}^l_A(A/A_{\ast,>d},-)$. As $A$ is tails coherent, it follows from Proposition 3.9(2) that $P_i/P_{i,>d} \in \gr A$ for all $i,d \in I$, and so the result follows from the previous proposition. 
\end{proof}

\begin{lemma}
If $N$ is a torsion module, then $\ssf{R}^l \tau(N) = 0$ for $l> 0$.
\end{lemma}
\begin{proof}
    Using the previous lemma we may suppose that $N$ is upper-bounded by some $u \in I$. Applying $\varinjlim_d\underline{\Hom}_A(-, N)$ to (3.4) gives us, for each $l\geq 0$, an exact sequence
    \[\varinjlim_d\underline{ \Ext}^l_A(A_{\ast,> d}, N) \to \ssf{R}^{l+1}\tau(N) \to 0.\]
    We show $\varinjlim_d \underline{\Ext}^l_A(A_{\ast,> d}, N) = 0$ for $l \geq 0$. For $i \in I$, $\varinjlim_d \underline{\Ext}^l_A(A_{\ast,> d}, N)_i = \varinjlim_d \Ext^l_A(P_{i,> d}, N)$. We now take a resolution $F^{\bullet}$ of $P_{i,> d}$ as in Proposition 4.1, where the positive indexing on $A$ ensures every summand $P_j$ of the free modules in $F^{\bullet}$ is such that $j > d$. Now $N$ is upper-bounded by $u$, so it follows that $\Ext^l_A(P_{i,> d}, N) = 0$ if $d \geq u$. Thus $\varinjlim_d\underline{\Ext}^l_A(A_{\ast,> d}, N)_i = 0$ for all $i \in I$ and all $l \geq 0$.
\end{proof}

In proving the previous lemma, the following result was also shown.

\begin{lemma}
    If $N$ is torsion, $\varinjlim_d\underline{\Ext}^l_A(A_{\ast,> d}, N)=0$ for all $l\geq 0$.
\end{lemma}

For $M \in \Gr A$, $\QHOM_A(\A,\MM)$ has the structure of a right $A$-module in the following way:
\begin{align*}
    \QHOM_A(\A,\MM)_i \times A_{ij} &= \QHom_A(\cP_i,\MM) \times \Hom_A(P_j,P_i)\\
    &\to \QHom_A(\cP_i,\MM) \times \QHom_A(\cP_j,\cP_i) \\
    &\to \QHom_A(\cP_j,\MM) \\
    &=\QHOM_A(\A,\MM)_j,
\end{align*}
where the first arrow is $id \times \pi$ and the second arrow is composition in $\QGr A$. In this way, $\QHOM_A(\A,-)$ becomes a functor $\QGr A \to \Gr A$. We also recall that given $M \in \Gr A$, the \textbf{saturation} of $M$ is the $A$-module $\widetilde{M} := \omega\pi M$. The proof of the next lemma is identical to the $\Z$-indexed case \cite[Lemma 3.17]{mori_char_2025}.

\begin{lemma}
    There is an isomorphism of functors $\omega(-) \cong \QHOM_A(\A,-):\QGr A \to \Gr A$, and hence $\ssf{R}^l\omega(\MM) \cong \underline{\E xt}^l_A(\A,\MM)$ in $\Gr A$ for all $M \in \Gr A$ and $l \geq 0$.
\end{lemma}

For the purposes of the next proof, we recall the construction of the morphisms in $\qgr A$. Given two coherent modules $M$ and $N$, we define a set
\[ J(M,N) := \{ (M',N') \mid M'\subset M,N'\subset N\text{ and } M/M',N' \in \ssf{tors}A\}.\]
This is a directed set, with $(M',N') < (M'',N'')$ if $M'' \subset M'$ and $N' \subset N''$, and the morphisms in $\qgr A$ are defined by the directed colimit
\[ \QHom_A(\MM,\NN) = \varinjlim_{(M',N') \in J(M,N)} \Hom_A(M',N/N'). \]

\begin{propn}
    For $M \in \Gr A$, $\widetilde{M} \cong \varinjlim_d\underline{\Hom}_A(A_{\ast,> d}, M)$ in $\Gr A$.
\end{propn}
\begin{proof}
    We first prove the case when $M \in \gr A$. We begin by showing that $\widetilde{M} \cong \varinjlim_d\underline{\Hom}_A(A_{\ast,> d}, M / \tau M)$. From Lemma 4.5 we have an isomorphism $\widetilde{M} \cong \QHOM_A(\A, \MM)$ in $\Gr A$. Then $\widetilde{M}_i \cong \QHom_A(\cP_i,\MM)$, which was described before the proposition. If $N \subset P_i$ and $P_i/N$ is upper-bounded by $d \in I$, $P_{i,> d} \subseteq N$. It follows that $\{ (P_{i,> d}, \tau M) ~|~ d \in I\}$ is cofinal in $J(P_i,M)$, and so $\widetilde{M}_i \cong \QHom_A(\cP_i, \MM) \cong \varinjlim_d\Hom_A(P_{i,> d}, M / \tau M)$. Thus we have right $A$-module isomorphisms
    \[ \widetilde{M} \cong \bigoplus_{i \in I} \varinjlim_d\Hom_A(P_{i,> d}, M / \tau M) \cong \varinjlim_d\underline{\Hom}_A(A_{\ast,> d}, M / \tau M).\]
    To conclude the proof we show $\varinjlim_d\underline{\Hom}_A(A_{\ast,> d}, M / \tau M) \cong \varinjlim_d\underline{\Hom}_A(A_{\ast,> d}, M)$. We apply $\varinjlim_d\underline{\Hom}_A(A_{\ast,> d}, - )$ to the exact sequence 
    \[0 \to \tau M \to M \to M / \tau M \to 0,\] 
    from which we obtain the desired isomorphism as Lemma 4.4 shows $\varinjlim_d\underline{\Ext}^l_A(A_{\ast,> d}, \tau M)=0$ when $l = 0, 1$. 
    
    As $A$ is tails-coherent, every $M \in \Gr A$ is a directed colimit of coherent modules $\varinjlim_j M^j$. The saturation functor commutes with directed colimits by Corollary 3.19(2) and the exactness of $\pi$. The modules $P_{i,>d}$ are coherent for all $d \in I$ as $A$ is tails-coherent, so $\varinjlim_d\underline{\Hom}_A(A_{\ast,> d}, -)$ also commutes with directed colimits. Thus we have right $A$-module isomorphisms
    \[ \widetilde{M} \cong \varinjlim_j \widetilde{M^j} \cong \varinjlim_j \varinjlim_d\underline{\Hom}_A(A_{\ast,> d}, M^j) \cong \varinjlim_d\underline{\Hom}_A(A_{\ast,> d}, M) \]
    for all $M \in \Gr A$.
\end{proof}

Applying $\varinjlim_d\HOM_A(-,M)$ to (3.4) then gives the following.

\begin{cor}
For every $M \in \mathsf{Gr}A$ there is an exact sequence 
\begin{align}
0 \to \tau M \to M \to \widetilde{M} \to \ssf{R}^1 \tau M \to 0,
\end{align}
and for every $l > 0$ there are isomorphisms $\ssf{R}^{l+1}\tau M \cong \ssf{R}^l\omega(\MM)$ in $\Gr A$.
\end{cor}

Recalling the notation of Corollary 3.19, we have a diagram of categories and functors:

\begin{center}
\begin{tikzcd}
\gr A \arrow[d, "\pi"] \arrow[r, hookrightarrow, "\iota"]
 & \mathsf{Gr}A \arrow[d, "\pi"', shift right] \\
\qgr A \arrow[r, hookrightarrow, "\underline{\iota}"]
 & \mathsf{QGr}A \arrow[u, "\omega"', shift right].
\end{tikzcd}
\end{center}
In many cases the section functor $\omega: \mathsf{QGr}A \to \mathsf{Gr}A$ does not restrict to a section functor $\omega:\qgr A \to \gr A$, and so while any coherent module $M$ has a saturation $\widetilde{M} = \omega \underline{\iota} \pi M$, the saturation of a coherent module need not be coherent.

\begin{defn}
For $M \in \gr A$ and $n \geq 0$, we say $M$ \textbf{satisfies} $\chi_n$ if for all $m \leq n$ and all $d \in I$, $\ssf{R}^{m} \tau(M)_{> d} \in \gr A$. If $M$ satisfies $\chi_n$ for all $M \in \gr A$ we say $A$ \textbf{satisfies} $\chi_n$, and if $\chi_n$ holds for all $n \geq 0$ we say $A$ \textbf{satisfies} $\chi$.
\end{defn}

The following characterisation of $\chi_1$ is similar to \cite[Proposition 3.14(2b)]{artin_noncommutative_1994}.

\begin{propn}
    For $M \in \gr A$, $M$ satisfies $\chi_1$ if and only if $\widetilde{M}_{> d}$ is coherent for all $d \in I$.
\end{propn}
\begin{proof}
    Taking the tails of (4.1) gives the exact sequence
    \[0 \to (\tau M)_{> d} \to M_{> d} \to \widetilde{M}_{> d} \to (\ssf{R}^1 \tau M)_{> d} \to 0.\]
    As $A$ is tails-coherent and $M \in \gr A$, $M_{>d}$ is also coherent by Proposition 3.9(3). Thus $\widetilde{M}_{> d}$ is coherent if and only if both $\tau(M)_{> d}$ and $\ssf{R}^1\tau(M)_{> d}$ are coherent, and so $\widetilde{M}_{> d}$ is coherent for all $d \in I$ if and only if $M$ satisfies $\chi_1$.
\end{proof}

\begin{remark}
    The definition given of the $\chi$ condition is based on that of \cite{artin_noncommutative_1994} for non-locally finite $\Z$-graded algebras. While our algebras are locally finite, we run into similar issues stemming from our index sets, as the tails $T_{> d}$ of an upper-bounded module $T \in \ssf{tors}A$ are not necessarily finite dimensional (c.f. Lemma 3.20). 
\end{remark}

\section{A characterisation of $\qgr$}
In this section we characterise abelian categories equivalent to $\qgr A$ for tails-coherent $I$-indexed algebras $A$, obtaining an $I$-indexed version of \cite[Theorem 2.4]{polishchuk_noncommutative_2005}. Throughout this section $\CC$ will be an abelian category with a sequence of objects $\E = (E_i)_{i \in I}$ indexed by $I$, and we insist $\End_{\CC}(E_i)$ is a division ring for each $i \in I$. Such a sequence defines a connected, positively $I$-indexed algebra 
\[A(\E) := \bigoplus_{i\leq j} \Hom_{\CC}(E_j,E_i).\] 
We will simply call this algebra $A$ when $\E$ is clear. There is a natural functor $\Gamma_{\ast}:\CC \to \Gr A$, which sends $X \in \CC$ to the right $A$-module $\bigoplus_{i \in I} \Hom_{\CC}(E_i, X)$. By definition we have that $\Gamma_{\ast}(E_i) = P_i \in \Gr A(\E)$ for all $i \in I$. We use $\Gamma_{> d}$ to denote the tails $(\Gamma_{\ast}(-))_{> d}$ for $d \in I$, and similarly for weak tails. 

\begin{defn}
Let $\E = (E_i)$ be an $I$-indexed sequence as above. We say $\E$ is
\begin{itemize}[noitemsep]
\item[(P)] \textbf{projective} if for each surjection $f:X \to Y$ in $\CC$, there is some $d_f \in I$ such that $\Hom_{\CC}(E_i, f)$ is surjective for $i > d_f$,
\item[(C)] \textbf{coherent} if $\Gamma_{> d}(X)$ is finitely generated for all $X \in \CC$ and all $d \in I$, 
\item[(A)] \textbf{ample} if $\E$ satisfies (P) and (C) and, for every $X \in \CC$ and every $d \in I$, there is a surjection $\bigoplus_{l=1}^r E_{j_l} \to X$ with each $j_l > d$.
\end{itemize}
\end{defn}

Recalling that $M_{>n} = M_{\geq n+1}$ for modules $M$ over a $\Z$-indexed algebra, \cite[Lemma 2.1]{polishchuk_noncommutative_2005} shows our notion of ample $I$-indexed sequences directly generalises ample $\Z$-indexed sequences as in \cite[Section 2]{polishchuk_noncommutative_2005}. 

We now fix $\E = (E_i)_{i \in I}$ to be an $I$-indexed sequence in an abelian category $\CC$.

\begin{propn}
    If $\E$ is coherent (C), then
    \begin{enumerate}
        \item $\Gamma_{>d}(X)$ is coherent for all $X \in \CC$ and $d \in I$. 
        \item $A(\E)$ is tails-coherent.
    \end{enumerate}  
\end{propn}
\begin{proof}
     (1): We know $\Gamma_{>d}(X)$ is finitely generated by (C), so we let $f: \bigoplus_{l=1}^r P_{j_l} \to \Gamma_{>d}(X)$ be any homomorphism. Note that each $j_l > d$ as $A$ is positively indexed. We then have 
    \[\Hom_A\left(\bigoplus_{l=1}^r P_{j_l},\Gamma_{>d}(X)\right) \cong \bigoplus_{l=1}^r \Gamma_{>d}(X)_{j_l} \cong \Hom_{\CC}\left(\bigoplus_{l=1}^r E_{j_l},X\right),\] 
    so there is some $g: \bigoplus_{l=1}^r E_{j_l} \to X$ in $\CC$ such that $\Gamma_{\ast}(g)= f$. As each $j_l > d$, $f = \Gamma_{> d}(g)$. Then, as $\Gamma_{> d}$ is left exact, we have
    \[\ker f = \ker \Gamma_{> d} (g) \cong \Gamma_{> d}(\ker g).\]
    Then $\Gamma_{> d}(\ker g)$ (and hence $\ker f$) are finitely generated by (C), so $\Gamma_{>d}(X)$ is coherent.
    
    (2) Coherence (C) shows $A$ is tails-finitely generated, as $P_{i,>d} = \Gamma_{>d}(E_i)$, and so it follows from Proposition 3.9(2) that $S_i = P_i/P_{i,>i}$ is finitely presented. Then $S_i$ is coherent, as any non-zero map $F \to S_i$ from a free and finitely generated module $F$ must be surjective, and the kernel of a surjective morphism between finitely presented modules is always finitely generated. Finally, $P_i$ is coherent as the exact sequence
    \[ 0 \to P_{i,>i} \to P_i \to S_i \to 0 \]
    expresses $P_i$ as an extension of the coherent module $S_i$ by the coherent module $P_{i,>i} = \Gamma_{>i}(E_i)$.
\end{proof}

It now follows from Theorem 3.18 and Corollary 3.19 that $\QGr A$ and $\qgr A$ exist when $\E$ is coherent. We will denote by $\QGamma_{\ast}$ the composition $\pi\Gamma_{\ast}$. We note that $\QGamma_{\ast}(X) \cong \pi \Gamma_{> d}(X)$ for any $X \in \CC$ and any $d \in I$, as the cokernel of the inclusion $\Gamma_{> d}(X) \hookrightarrow \Gamma_{\ast}(X)$ is upper-bounded by $d$. Each tail $\Gamma_{> d}(X)$ is coherent via the previous proposition, so $\QGamma_{\ast}:\CC \to \qgr A$. 

\begin{lemma}
    If $\E$ is projective (P) and coherent (C), $\QGamma_{\ast}$ is exact.
\end{lemma}
\begin{proof}
    Coherence (C) tells us that $\QGr A$ and $\qgr A$ exist. As $\Gamma_{\ast}$ is left exact and $\pi$ is exact, $\QGamma_{\ast}$ is left exact too. Given a surjection $g$ in $\CC$, (P) equivalently says that the cokernel of $\Gamma_{\ast}(g)$ is upper-bounded. Then $\coker\QGamma_{\ast}(g) \cong \pi (\coker \Gamma_{\ast}(g)) \cong 0$, so $\QGamma_{\ast}$ is also right exact.
\end{proof}

\begin{lemma}
    If $\E$ is ample (A), $\Gamma_{\ast}(X)$ is torsion-free for all $X \in \CC$.
\end{lemma}
\begin{proof}
    Suppose $m \in \Gamma_{i}(X)$ is torsion, with $(mA)_{> d} = 0$ for some $d \in I$. Using (A) we obtain a surjection $g_d: \bigoplus_{l=1}^r E_{j_l} \to E_i$ with each $j_l > d$, so that $g_d \in P_{i, > d}$. Then $m g_d = 0$, and as $g_d$ is surjective we must have $m = 0$.
\end{proof}

\begin{propn}
    If $\E$ is ample (A), $\Gamma_{\ast}(X) \cong \widetilde{\Gamma_{\ast}(X)}$ for all $X \in \CC$.
\end{propn}
\begin{proof}
    The previous lemma and (4.1) give us an exact sequence
    \[0 \to \Gamma_{\ast}(X) \to \widetilde{\Gamma_{\ast}(X)} \to \ssf{R}^1\tau \Gamma_{\ast}(X) \to 0.\]
    We show $\ssf{R}^1\tau\Gamma_{\ast}(X) = 0$. Suppose $\overline{m} \in \ssf{R}^1\tau \Gamma_{\ast}(X)_d$ is the image of $m \in \widetilde{\Gamma_{\ast}(X)}_d$. As $\ssf{R}^1\tau\Gamma_{\ast}(X)$ is torsion, there is some $u \in I$ for which $(mA)_{> u} \subseteq \Gamma_{\ast}(X)$. We now use (A) to construct a presentation of $E_d$
    \[\bigoplus_{t=1}^s E_{i_t} \xrightarrow{f} \bigoplus_{l=1}^{r} E_{j_{l}} \xrightarrow{g} E_d \to 0\]
    with all $j_l,i_t > u$. We write these maps as $f=(f_{lt})$ and $g=(g_{l})$. We also assume that $(\sum_{l=1}^r g_l A)_{>d'} = P_{d,>d'}$ for some $d'>d$. This is possible as $A$ is tails-coherent by Proposition 5.2(2), so $P_{d,>d'}$ is finitely generated, and thus we may simply add a finite set of generators for $P_{d,>d'}$ to $\{g_l\}$. Applying $\Hom_{\CC}(-,X)$ gives the exact sequence
    \[ 0 \to \Hom_{\CC}(E_d, X) \xrightarrow{\Hom_{\CC}(g,X)} \bigoplus_{l=1}^r\Hom_{\CC}(E_{j_l},X) \xrightarrow{\Hom_{\CC}(f,X)} \bigoplus_{t=1}^s\Hom_{\CC}(E_{i_t},X). \]
    We have $mg_l \in \Gamma_{\ast}(X)$ for $1 \leq l \leq r$ as each $j_l>u$, and hence $mg \in \bigoplus_{l=1}^r \Hom_{\CC}(E_{j_l},X)$. Note that this does not immediately imply that $m \in \Gamma_{\ast}(X)$, however, as for $j \not>u$ there may be $x \in A_{dj}$ with $mx \not\in \Gamma_{j}(X)$. From the exact sequence above we see that $(mg)f = m(gf) = 0 \in \bigoplus_{t=1}^s\Hom_{\CC}(E_{i_t},X)$, and so it follows that there is some $m' \in \Hom_{\CC}(E_d,X) = \Gamma_{d}(X)$ such that $m'g = mg$. Then $\sum_{l=1}^r(m'-m)g_l=0$, and as we chose $\{g_l\}$ to contain a set of generators for $P_{d,>d'}$ we see that $((m'-m)A)_{>d'} = 0$. So the element $m'-m \in \widetilde{\Gamma_{\ast}(X)}$ is torsion, but $\widetilde{\Gamma_{\ast}(X)}$ is torsion-free by definition, so we must have $m'=m$. Thus $\overline{m} = 0$, and hence $\ssf{R}^1\tau\Gamma_{\ast}(X)=0$.
\end{proof}

\begin{cor}
    If $\E$ is ample (A), $\QHom_A(\QGamma_{\ast}(X),\QGamma_{\ast}(Y)) \cong \Hom_A(\Gamma_{\ast}(X),\Gamma_{\ast}(Y))$ for all $X,Y \in \CC$.
\end{cor}
\begin{proof}
    This follows from the adjoint isomorphism and the previous proposition, as
    \[ \QHom_A(\QGamma_{\ast}(X),\QGamma_{\ast}(Y)) \cong \Hom_A(\Gamma_{\ast}(X),\widetilde{\Gamma_{\ast}(Y)}) \cong \Hom_A(\Gamma_{\ast}(X),\Gamma_{\ast}(Y)). \]
\end{proof}

\begin{propn}
    If $\E$ is ample (A), $\QGamma_{\ast}$ is fully faithful.
\end{propn}
\begin{proof}
    If $X = E_i$ we have
    \[ \Hom_{\CC}(X,Y) = \Gamma_{i}(Y) \cong \Hom_A(P_i,\Gamma_{\ast}(Y)). \]
    It follows from the previous corollary that $\Hom_{\CC}(E_i,Y) \cong \QHom_A(\cP_i,\QGamma_{\ast}(Y))$. The same then must be true when $X$ is a finite sum of the $E_i$. The general result follows by applying $\Hom_{\CC}(-,Y)$ to a resolution of $X$ by finite sums of the $E_i$, which exists via (A).
\end{proof}

\begin{propn}
    If $\E$ is ample (A), $\QGamma_{\ast}$ is essentially surjective.
\end{propn}
\begin{proof}
    For any $\MM = \pi M \in \qgr A$, we take a finite presentation of $M$
    \[\bigoplus_{t=1}^s P_{i_t} \xrightarrow{f} \bigoplus_{l=1}^r P_{j_l} \to M \to 0.\]
    Now $\QGamma_{\ast}$ is fully faithful by the previous proposition, so there is a morphism $g: \bigoplus_{t=1}^s E_{i_t} \to \bigoplus_{l=1}^r E_{j_l}$ such that $\QGamma_{\ast}(g) \cong \pi f$. Lemma 5.3 shows $\QGamma_{\ast}$ is exact, so 
    \[\QGamma_{\ast}(\coker g) \cong \coker \QGamma_{\ast}(g) \cong \coker \pi f \cong \pi (\coker f) \cong \MM,\]
    and hence $\QGamma_{\ast}$ is essentially surjective.
\end{proof}

\begin{propn}
    If $\E$ is ample (A), $A(\E)$ satisfies $\chi_1$.
\end{propn}
\begin{proof}
    Firstly, $A$ is tails-coherent by Proposition 5.2(2). For any $M \in \gr A$, Propositions 5.5 and 5.8 show there is some $X \in \CC$ with $\Gamma_{\ast}(X) \cong \widetilde{M}$. The tails $\Gamma_{> d}(X) \cong \widetilde{M}_{> d}$ are coherent for all $d \in I$ by Proposition 5.2(1), so $M$ satisfies $\chi_1$ via Proposition 4.9.
\end{proof}

We are now able to prove an $I$-indexed version of Polishchuk's theorem \cite[Theorem 2.4]{polishchuk_noncommutative_2005}.

\begin{theorem}
    If $\CC$ has an ample $I$-indexed sequence $\E = (E_i)_{i \in I}$, then $A(\E)$ is tails-coherent, satisfies $\chi_1$, and $\QGamma_{\ast}:\CC \to \qgr A(\E)$ is an exact equivalence. Conversely, $(\cP_i)_{i \in I}$ is an ample $I$-indexed sequence in $\qgr A$ for any tails-coherent $I$-indexed algebra $A$ satisfying $\chi_1$.
\end{theorem}
\begin{proof}
    If $\E$ is an ample sequence, the results of this section show $A(\E)$ is tails-coherent, satisfies $\chi_1$, and that $\QGamma_{\ast}$ defines an exact equivalence $\CC \simeq \qgr A$. 
    
    For the converse, suppose $A$ is a tails-coherent $I$-indexed algebra satisfying $\chi_1$ and consider the $I$-indexed sequence $(\mathcal{P}_i)_{i \in I}$. In this case $\Gamma_{\ast} = \omega(-)$ and $\QGamma_{\ast} = \pi\omega$ is naturally isomorphic to the identity functor. We fix $\MM = \pi M$ and $\NN = \pi N$ for some $M,N \in \gr A$.

    (P): Let $f: \MM \to \NN$ be a surjection in $\qgr A$, and let $C$ denote the cokernel of $\Gamma_{\ast}(f)=\omega(f):\widetilde{M} \to \widetilde{N}$. We show $C$ is upper-bounded by some $u \in I$, so that $\QHom_A(\cP_i,f)$ is surjective for $i\geq u$. Note that $C$ is torsion, as $\pi C \cong \coker f = 0$. As $A$ satisfies $\chi_1$, both $\widetilde{M}_{>d}$ and $\widetilde{N}_{>d}$ are coherent for any $d \in I$ by Proposition 4.9, so $C_{>d} = \coker \Gamma_{>d}(f)$ is coherent. It follows that $C_{>d}$ (and hence also $C$) is upper-bounded, so $(\cP_i)_{i \in I}$ is projective (P).
    
    (C): As $A$ satisfies $\chi_1$, Proposition 4.9 ensures $\Gamma_{> d}(\MM) \cong \widetilde{M}_{> d}$ is coherent for any $d \in I$. Hence $(\cP_i)_{i \in I}$ is coherent (C). 
    
    (A): Take any $d \in I$. The tail $M_{> d}$ of $M$ is finitely generated by Proposition 3.9(3), so there is a surjection
    \[ \bigoplus_{l=1}^r P_{j_l} \to M_{> d}. \]
    Note that each $j_l>d$, as $A$ is positively indexed. Applying $\pi$ gives a surjection
    \[\bigoplus_{l=1}^r \cP_{j_l} \to \pi (M_{>d}) \cong \MM, \]
    and so $(\cP_i)_{i \in I}$ is ample (A).
\end{proof}

\section{Examples of tails-finitely generated $I$-indexed algebras}
We now look at some of the specific examples mentioned in the introduction. We will see that the relevant graded/indexed algebras in these examples are tails-finitely generated, and that the associated noncommutative projective schemes constructed for each grading/indexing set coincide with the constructions for $I$-indexed algebras from Section 3. In this way, all of these examples can be treated uniformly using the theory of $I$-indexed algebras developed in the earlier sections. In particular, the noncommutative projective schemes in each case can be characterised in terms of ample $I$-indexed sequences as in Theorem 5.10, providing a method to identify/construct such noncommutative projective schemes which is applicable to all of these index sets.

\subsection{Twisted multi-homogeneous coordinate rings}
One large class of rings arising in noncommutative algebraic geometry is the \textit{twisted homogeneous coordinate rings}, introduced by Artin and Van den Bergh in \cite{artin_twisted_1990}. Multi-homogeneous versions of these algebras, constructed by Chan in \cite{chan_twisted_2000}, are $\N^r$-graded rings of the form 
\[R=\bigoplus_{(m)=(m_1,\ldots,m_r) \in \N^r} \ssf{H}^0(X;\L_1^{m_1} \otimes \L_2^{m_2} \otimes \ldots \otimes \L_r^{m_r}),\]
where $X$ is a projective scheme over an algebraically closed field $k$, and $\L_1,\ldots,\L_r$ are invertible $\O_X$-bimodules satisfying some commutativity relations and an ampleness condition \cite[Definitions 3.1,~3.2]{chan_twisted_2000}. Such a ring is $\N^r$-graded, and so the associated $\Z^r$-indexed algebra $B = \mathbf{B}_{\Z^r}R$ is positively $\Z^r$-indexed. 

To see that $B$ is tails-finitely generated we use Proposition 3.12(2a). We let $1_i\in \Z^r$ denote the element with 1 in the $i$-th position and 0 in all others. For any $(m)\in\Z^r$, $e_{(m)}\m =P_{(m),>(m)}$ is finitely generated `in degrees $1_i$' by $B_{(m),(m)+1_i} = \ssf{H}^0(X;\L_i)$ for $i=1,\ldots,r$, so $B$ satisfies $(\ast)$. The commutation relations give $k$-bimodule isomorphisms 
\[B_{(m),(m)+1_i}B_{(m)+1_i,(m)+1_i+1_j} \cong B_{(m),(m)+1_j}B_{(m)+1_j,(m)+1_i+1_i}\]
for all $(m)\in \Z^r$ and $1 \leq i,j \leq r$. Thus $B_{(m)(m')}B_{(m')(n)} \cong B_{(m)(n)}$ for all $(m)\leq (m') \leq (n)$, and so the multiplication maps $\mu_{(m)(m')(n)}: B_{(m)(m')} \otimes_k B_{(m')(n)} \to B_{(m)(n)}$ are all surjective. In Example 3.14(1) we saw any $(m)$ and $(n)$ in $\Z^r$ have a least upper bound so, in the notation of Proposition 3.12, $\vert U_{min}((m),(n)) \vert =1$. It then follows from Proposition 3.12(2a) that (the $\Z^r$-indexed algebra associated to) a twisted multi-homogeneous coordinate ring must be tails-finitely generated. 
    
The subcategory of graded torsion modules \cite[Definition 3.3]{chan_twisted_2000} for twisted multi-homogeneous coordinate rings is also defined in terms of directed colimits of upper-bounded modules (referred to as \textit{first quadrant bounded} in \cite{chan_twisted_2000}), so the equivalence $\Gr R \simeq \Gr B$ given by \cite[Theorem 3.6]{sierra_g-algebras_2011} shows the two notions of torsion agree. It then follows from \cite[Theorem 3.4]{chan_twisted_2000} that we have an equivalence
\[ \QGr B \simeq \ssf{Qcoh}X. \]
Thus if $R$ is coherent (a criterion for $R$ to be noetherian is given in \cite[Theorem 5.2]{chan_twisted_2000}), $\qgr B$ coincides with the category $\coh X$ of coherent sheaves on $X$.

\subsection{$\Z^2$-indexed algebras and noncommutative birational transformations}
Another example of $\Z^2$-indexed algebras arising in noncommutative algebraic geometry comes from the study of noncommutative versions of the classical birational transformation $\PP^1(k) \times \PP^1(k) \dashrightarrow \PP^2(k)$ \cites{presotto_birational_2018, presotto__2017}. The $\Z^2$-indexed algebras in question (which we will call $A$) are constructed by `gluing' the $\Z$-indexed algebras $\mathbf{B}_{\Z}R$ and $\mathbf{B}_{\Z}S$ associated to a quadratic Sklyanin algebra $R$ and a cubic Sklyanin algebra $S$ \cite[Sections 7.1,~7.2]{presotto_birational_2018}. We view these algebras as the homogeneous coordinate rings of a noncommutative $\PP^2$: $\X=\qgr (\mathbf{B}_{\Z}R)$ and a noncommutative $\PP^1 \times \PP^1$: $\Y=\qgr (\mathbf{B}_{\Z}S)$. The corresponding glued $\Z^2$-indexed algebra $A$ is related to a common noncommutative ($\Z$-indexed) blow-up $\mathbb{W}=\qgr T$ of $\X$ and $\Y$ by an equivalence $\mathbb{W} \simeq \qgr A$ \cite[Theorem 1.1]{presotto__2017}. 

We again use Proposition 3.12 to see that the glued $\Z^2$-indexed algebra $A$ is tails-finitely generated. In this case we use Proposition 3.12(1), as $A$ is always noetherian \cite[Theorem 5.5]{presotto__2017}. We also note that the definition of torsion for $\Z^2$-indexed algebras as in \cite[Definition 3.3]{presotto__2017} agrees with that of Definition 3.16 for $I=(\Z^2,\leq_{prod})$, so again the associated noncommutative projective schemes coincide.

\subsection{Weighted projective lines}
Leaving the realm of $\N^r$-graded rings, we look at the weighted projective lines introduced by Geigle and Lenzing \cite{geigle_class_1987}, which are graded by a free abelian group of rank 1. We fix an algebraically closed field $k$, a finite sequence of distinct points $\underline{p}=(p_0,\ldots,p_n)$ in $\PP^1_k$, and a sequence of weights $\underline{w}=(w_0,\ldots,w_n) \in (\Z^+)^n$. Without loss of generality we may suppose $p_0=\infty$, $p_1=0$ and $p_2=1$. 

We let $L(\underline{w})$ be the rank 1 abelian group with generators $x_0,\ldots,x_n$ and relations $w_0x_0=\ldots=w_nx_n$, and let $S(\underline{w})$ denote the $L(\underline{w})$-graded polynomial ring $k[X_0,\ldots,X_n]$ with $\deg X_i=x_i$. The proj $\mathsf{Proj}(S(\underline{w}))$ is a \textit{weighted projective space} $\PP^n(\underline{w})$ \cite[Section 1]{geigle_class_1987}. From the sequence of points $\underline{p}$ we define a homogeneous ideal $I(\underline{w},\underline{p})$ of $S(\underline{w})$ generated by the elements $X_i^{w_i}-(X_1^{w_1}-p_iX_0^{w_0})$ for $i=2,\ldots,n$. The \textit{weighted projective line} $C(\underline{w},\underline{p})$ is then defined as the proj $\mathsf{Proj}(S(\underline{w},\underline{p}))$ of the quotient ring $S(\underline{w},\underline{p})=S(\underline{w})/I(\underline{w},\underline{p})$, and is naturally viewed as a curve in the weighted projective space $\PP^n(\underline{w})$.

As we have seen in Examples 3.1(2) and 3.14(2), $L(\underline{w})$ has a natural partial ordering under which it is locally finite and directed. Each generator $X_i$ for $S(\underline{w},\underline{p})$ satisfies $\deg(X_i)=x_i>0$, so $S(\underline{w},\underline{p})$ is positively $L(\underline{w})$-graded. Thus $D:= \mathbf{B}_{L(\underline{w})}S(\underline{w},\underline{p})$ is a positively $L(\underline{w})$-indexed algebra. Moreover, $D$ is noetherian as $S(\underline{w},\underline{p})$ is, and hence is tails-coherent by Proposition 3.12(1).

Theorem 1.8 of \cite{geigle_class_1987} gives the following version of Serre's theorem:
\[ \coh(C(\underline{w},\underline{p})) \simeq \gr (S(\underline{w},\underline{p}))/\gr_0 (S(\underline{w},\underline{p})), \]
where $\gr_0 (S(\underline{w},\underline{p}))$ is the Serre subcategory of finite length $L(\underline{w})$-graded $S(\underline{w},\underline{p})$-modules. Clearly $\ssf{fd}D$ coincides with $\gr_0 (S(\underline{w},\underline{p}))$. As the partial order on $L(\underline{w})$ satisfies the condition of Proposition 3.12(2b), Lemma 3.20 shows $\ssf{tors}D = \ssf{fd} D$, and so we have equivalences
\[ \qgr D \simeq \gr (S(\underline{w},\underline{p}))/\gr_0 (S(\underline{w},\underline{p})) \simeq \coh (C(\underline{w},\underline{p})).\]

\section{Preprojective algebras of non-algebraic species}
In this section we provide an application of Theorem 5.10 to the representation theory of hereditary artinian rings. Representations of bimodule species generalise the representation theory of quivers, and in a similar way capture the representation theory of a large class of hereditary artinian rings (including hereditary finite dimensional algebras over perfect fields) \cite{dlab_indecomposable_1976}. A species is \textit{algebraic} (or a $k$\textit{-species}) if its path algebra is a finite dimensional algebra over a field $k$; otherwise the species is \textit{non-algebraic}. Connections between the representation theory of species and coherent sheaves on noncommutative projective schemes are well known for algebraic species \cite{baer_tame_1987, lenzing_curve_1986, minamoto_ampleness_2012,krause_serre_2024}. Less is known for non-algebraic species, with the exception of species over the Kronecker quivers \cite{chan_species_2016,chan_representation_2019}. We use Theorem 5.10 to construct a noncommutative projective scheme derived-equivalent to the path algebra of any bimodule species of infinite representation type. Moreover, we will see that the preprojective algebra of the species may be recovered as a `noncommutative Veronese' of the indexed algebra we construct, similarly to \cite[Proposition 7.5]{chan_species_2016}. 

Let $Q=(Q_0,Q_1,d)$ be a valued quiver, which we will insist is finite, connected, acyclic and non-Dynkin. A \textbf{modulation} $\mathfrak{M}=(\{K_v\}_{v \in Q_0},\{B_{e}\}_{e \in Q_1})$ of $Q$ is an assignment of division rings $K_v$ for each $v \in Q_0$ and $K_v-K_w$-bimodules $B_e$ for every edge $e:v \to w \in Q_1$ satisfying \begin{enumerate}
    \item[(M1)] $\dim_{K_v} B_e = d^e_{wv}$ and $\dim_{K_w} B_e = d^e_{vw}$,
    \item[(M2)] ${}^{\ast}B_e := \Hom_{K_v}(B_e,K_v) \cong \Hom_{K_w}(B_e,K_w)=:B_e^{\ast}$ as $K_w-K_v$-bimodules.
\end{enumerate}
A \textbf{bimodule species} $\cS$ is a pair $(Q,\mathfrak{M})$ consisting of a valued quiver $Q$ and a modulation $\mathfrak{M}$ of $Q$. In practice we may suppose $Q$ is without multiple edges by considering the bimodules $B_{vw} := \bigoplus_{e:v \to w \in Q_1} B_e$. Note that from the condition (M2) we have $B_{vw}^{\ast\ast} \cong {}^{\ast}(B_{vw}^{\ast}) \cong B_{vw}$ for all $v,w \in Q_0$.

A species $\cS$ has an associated `path algebra' $\Lambda=\Lambda(\cS)$ defined as the tensor algebra $T_K(B)$, where $K=\prod_{v \in Q_0}K_v$ and $B$ is the $K$-bimodule $\bigoplus_{v,w \in Q_0} B_{vw}$. We will denote the category of finitely presented right $\Lambda$-modules by $\mod \Lambda$, and will write $D^b(\Lambda)$ to mean $D^b(\mod \Lambda)$. Similarly, $\Lambda\mod$ will denote the category of finitely presented left $\Lambda$-modules. For $v \in Q_0$, the unit of the division ring $K_v$ determines an idempotent $\varepsilon_v$ in $\Lambda$, and we denote by $P_v$ the indecomposable projective $\Lambda$-module $\varepsilon_v\Lambda$. The simple $\Lambda$-modules have the form $\varepsilon_v \Lambda \varepsilon_v$, and will be denoted by $S_v$. 

The ring $\Lambda$ has a duality $D(-):\mod\Lambda \to \Lambda\mod$ represented by the dual bimodule $D\Lambda := \Hom_{K}(\Lambda_K,K)$. This is naturally an $\End_{\Lambda}(D\Lambda_{\Lambda})-\Lambda$-bimodule, and we identify $\End_{\Lambda}(D\Lambda_{\Lambda})$ with $\Lambda$ via the isomorphisms
\[ \Hom_{\Lambda}(D\Lambda_{\Lambda},D\Lambda_{\Lambda}) \cong \Hom_K(D\Lambda \otimes_{\Lambda} \Lambda_K,K) \cong \Lambda. \]
Here, the final isomorphism is determined by the double dual isomorphisms obtained from (M2). The functor $\Hom_{\Lambda}(-,D\Lambda_{\Lambda})\simeq\Hom_K(-\otimes_{\Lambda}\Lambda_K,K)$ is exact, so $D\Lambda$ is injective. Moreover, the indecomposable injective $\Lambda$-modules are $E_v:=\varepsilon_vD\Lambda$ for $v \in Q_0$. The simple modules $S_v \cong K_v$ are self-dual, so dualising the projective cover ${}_{\Lambda}(\Lambda\varepsilon_v) \to {}_{\Lambda}S_v$ gives an injective envelope $D({}_{\Lambda}S_v)=(S_v)_{\Lambda} \hookrightarrow (\varepsilon_vD\Lambda)_{\Lambda}$. Then $\sigma:=D\Lambda[-1]\in D^b(\Lambda)$ is a \textit{canonical complex} in the sense of \cite[Definition 3.0.1]{chan_representation_2019}, and hence $- \dtensor_{\Lambda} \sigma$ is an autoequivalence of $D^b(\Lambda)$ with inverse $-\dtensor_{\Lambda} \sigma^{-1} = \RHom_{\Lambda}(\sigma,-)$. We will denote the action of $\sigma$ on $X \in D^b(\Lambda)$ by $\sigma X$.

The ring $\Lambda$ is always a basic, hereditary artinian ring such that $\mod \Lambda$ is \textit{Krull-Schmidt} (any $M \in \mod \Lambda$ is isomorphic to a finite direct sum of indecomposable modules). Any such ring has 3 types of indecomposable modules: the \textit{preprojective}, \textit{postinjective} and \textit{regular} modules \cite{dlab_indecomposable_1976}. These modules are usually defined in terms of the \textit{Auslander-Reiten translate} $\tau = D \Ext_{\Lambda}^1(-,\Lambda)$ and its left adjoint $\tau^- = \Ext^1_{\Lambda}(D\Lambda,-)$. We recognise $\tau^{-}$ as $H^0(-\dtensor_{\Lambda}\sigma^{-1})$, as
\[ H^0(-\dtensor_{\Lambda}\sigma^{-1})=H^0\RHom_{\Lambda}(\sigma,-) = \Hom_{D^b(\Lambda)}(D\Lambda[-1],-) \cong \Ext^1_{\Lambda}(D\Lambda,-),\]
and so we may define these types of modules in terms of $\sigma$ as follows.

\begin{defn}
    A module $M \in \mod\Lambda$ is \begin{enumerate}
        \item \textbf{preprojective} if $H^0(\sigma^{n}M) = 0$ for $n \gg 0$.
        \item \textbf{postinjective} if $H^0(\sigma^{-n}M)=0$ for $n \gg 0$.
        \item \textbf{regular} if $\sigma^n M \in \mod\Lambda$ for all $n \in \Z$.
    \end{enumerate}
\end{defn}

We denote the full subcategories of preprojective, postinjective, and regular modules by $\ssf{P},\ssf{I}$ and $\ssf{R}$, respectively. We recall the following facts about $\mod \Lambda$: \begin{itemize}
    \item $\Hom_{\Lambda}(\ssf{I},\ssf{P})=\Hom_{\Lambda}(\ssf{I},\ssf{R})=\Hom_{\Lambda}(\ssf{R},\ssf{P})=\{0\}$.
    \item For any $M \in \ssf{P}$, there are finitely many indecomposable $N \in \mod \Lambda$ with $\Hom_{\Lambda}(N,M) \neq 0$.
    \item For any $M \in \ssf{I}$, there are finitely many indecomposable $N \in \mod \Lambda$ with $\Ext^1_{\Lambda}(N,M) \neq 0$.
\end{itemize}
These three subcategories determine a decomposition $\mod \Lambda = \ssf{P} \vee \ssf{R} \vee \ssf{I}$, where the notation $\CC_1 \vee \CC_2$ for disjoint subcategories $\CC_1,\CC_2$ of a category $\CC$ means the full subcategory consisting of direct sums $C_1 \oplus C_2$ with $C_1 \in \CC_1,~C_2\in \CC_2$. As $\Lambda$ is hereditary, the decomposition of $\mod \Lambda$ induces a decomposition 
\[D^b(\Lambda) = \bigvee_{n \in \Z} \ssf{P}[n] \vee \ssf{R}[n] \vee \ssf{I}[n],\]
as any $X \in D^b(\Lambda)$ may be identified with $\bigoplus_{n \in \Z}H^n(X)[-n]$. We denote the standard $t$-structure on $D^b(\Lambda)$ by $(\DD^{\leq 0},\DD^{\geq 0})$. As $\Hom_{\Lambda}(\ssf{I},\ssf{P}\vee \ssf{R})=0$ and $\mod \Lambda = \ssf{P}\vee \ssf{R} \vee \ssf{I}$, \cite[Theorem 2.1]{dickson_torsion_1966} shows $(\ssf{I},\ssf{P}\vee \ssf{R})$ is a torsion pair in $\mod \Lambda$. This torsion pair induces a $t$-structure $(\CC^{\leq 0},\CC^{\geq 0})$ on $D^b(\Lambda)$ with heart $\H = \ssf{I}[-1]\vee \ssf{P} \vee \ssf{R}$, and for any $X,Y \in \H$ we have isomorphisms $\Hom_{D^b(\H)}(X,Y[l]) \cong \Hom_{D^b(\Lambda)}(X,Y[l])$ when $l \in \{0,1\}$ \cite[Corollary 2.2]{happel_tilting_1996}. Furthermore, $D^b(\Lambda)$ and $D^b(\H)$ are triangle-equivalent \cite[Theorem 3.3]{happel_tilting_1996}. We will frequently use these facts without explicit mention. The $t$-structure defining $\H$ is given by the sets 
\[\CC^{\leq 0} = \{ X \in \DD^{\leq 1} \mid H^1(X) \in \ssf{I}\},\qquad \CC^{\geq 0} = \{ X \in \DD^{\geq 0} \mid H^{0}(X) \in \ssf{P} \vee \ssf{R}\}.\]

\begin{lemma}
    \begin{enumerate}
        \item The restriction $\sigma|_{\H}$ is an autoequivalence.
        \item The set of indecomposable objects in $\ssf{I}[-1]$ (resp. $\ssf{P}$) is $\{\sigma^mP_v\mid m > 0,~v \in Q_0\}$ (resp. $\{\sigma^mP_v\mid m \leq 0,~v \in Q_0\}$).
        \item For any $X \in \H$, $\sigma^{-n}X \in \ssf{P} \vee \ssf{R} \subset \mod \Lambda$ for $n \gg 0$.
    \end{enumerate}
\end{lemma}
\begin{proof}
    Claims (1) and (2) follow directly from \cite[Lemma 3.3(1) and Lemma 3.3(2)]{chan_representation_2019}. If we write $X$ as a sum $\bigoplus_{l=1}^r X_r$ of indecomposable objects, then from (2) we see that each summand $X_l$ is either regular or of the form $\sigma^{m_l}P_{v_l}$ for some $(m_l,v_l) \in \Z \times Q_0$. By definition we have that $\sigma^{m}(\ssf{R}) \subset \ssf{R}$ for all $m \in \Z$ and $\sigma^{-m}P_v \in \ssf{P}$ for all $m \geq 0$ and all $v \in Q_0$, so $\sigma^{-n} X \in \ssf{P}\vee \ssf{R}$ for sufficiently large $n$.
\end{proof}

We let $\E$ be the $\Z Q$-indexed sequence $(\sigma^mP_v)_{(m,v) \in \Z Q}$ in $\H$, where $\Z Q = (\Z \times Q_0,\prec)$ is the locally finite directed poset defined in Examples 3.1(3) and 3.14(3), and let $A = A(\E)$. As $\sigma$ is an autoequivalence we have
\begin{equation}
    A_{(m,v)(n,w)} = \Hom_{\H}(\sigma^nP_w,\sigma^mP_v) \cong \Hom_{\H}(\sigma^{n+1}P_w,\sigma^{m+1}P_v) = A_{(m+1,v)(n+1,w)}
\end{equation}
for all $(m,v),(n,w) \in \Z Q$. We also mention that $A$ satisfies $(\ast)$. It suffices to see that $P_{(0,v),>(0,v)}$ is finitely generated for each $v \in Q_0$ via (7.1). For each $v \in Q_0$, $P_{(0,v),>(0,v)}$ is finitely generated by any finite right basis of $\bigoplus_{w \in Q_0} B_{vw} \oplus B_{wv}^{\ast}$, and such a basis exists by the condition (M1) for the modulation $\mathfrak{M}$.

\begin{theorem}
    Let $\cS$ be a bimodule species with path algebra $\Lambda$, and let $\E$ and $A$ be as above. Then \begin{enumerate}
        \item $A$ is tails-coherent, so $\PP(\cS):=\qgr~A$ exists.
        \item $D^b(\Lambda) \simeq D^b(\PP(\cS))$.
        \item $\PP(\cS)$ is a hereditary category.
    \end{enumerate}
\end{theorem}
\begin{proof}
    We show that $\E$ is an ample $\Z Q$-indexed sequence, from which (1) will follow by way of Theorem 5.10. We have an equivalence $D^b(\H) \simeq D^b(\Lambda)$ by construction, and $\H$ is hereditary as $\Lambda$ is, so an exact equivalence $\H \simeq \PP(\cS)$ would induce an equivalence $D^b(\PP(\cS)) \simeq D^b(\Lambda)$ and show that $\PP(\cS)$ is hereditary. Thus (1)-(3) all follow from Theorem 5.10 if $\E$ is an ample sequence. 
    
    We recall that a surjection $f:X \to Y \in \H$ is a morphism $f:X \to Y \in D^b(\Lambda)$ with cone $C\in \CC^{\leq -1}=\DD^{\leq -1} \cup \ssf{I}$. 
    
    (P): Take a surjection $f:X \to Y \in \H$, and let $C$ be a cone for $f$ in $D^b(\Lambda)$ with $C \in \CC^{\leq -1}$. Without loss of generality we suppose $C$ is indecomposable. If $C \in \DD^{\leq -1}$ then $\Gamma_{\ast}(C)=0$, so we suppose $C \in \ssf{I}$. Then $\Gamma_{(m,v)}(C) \in \Hom_{\H}(\ssf{I}[-1],C) \cong \Ext^1_{\Lambda}(\ssf{I},C)$ for $m>0$. As $C \in \ssf{I}$, there are only finitely many $N \in \ssf{I}$ with $\Ext^1_{\Lambda}(N,C) \neq 0$. Then there are only finitely many $(m,v)$ with $m>0$ such that $\Gamma_{(m,v)}(C) \neq 0$, and so $\Gamma_{\ast}(C)$ is upper bounded. Thus $\E$ is projective (P).

    (C): As $A$ satisfies $(\ast)$, it follows from Remark 3.15 and Proposition 3.12(2b) that $A$ is tails-finitely generated. Thus $P_{(m,v),>(n,w)}=\Gamma_{>(n,w)}(\sigma^mP_v)$ is finitely generated for all $(m,v),(n,w) \in \Z Q$, and hence $\Gamma_{>(n,w)}(X)$ is finitely generated for all $X \in \ssf{I}[-1] \vee \ssf{P}$. It remains to see that $\Gamma_{>(n,w)}(X)$ is finitely generated for $X \in \ssf{R}$. As $X \in \ssf{R}$ is finitely generated we may take a finitely generated projective $\Lambda$-module $P$ and a surjection $f:P \to X$. As both $P$ and $X$ also belong to $\H$, $f$ is a surjection in $\H$ too. Consider the exact sequence 
    \[ 0 \to K \to \Gamma_{>(n,w)}(X) \to C \to 0, \]
    where $K = \mathrm{im}~\Gamma_{>(n,w)}(f)$ and $C=\coker \Gamma_{>(n,w)}(f)$. Now $K$ is finitely generated as it is a quotient of the finitely generated $A$-module $\Gamma_{>(n,w)}(P)$, so we see $\Gamma_{>(n,w)}(X)$ is finitely generated if and only if $C$ is. 
    
    Note that $\Gamma_{\ast}(X)$ is locally finite, as $\sigma^m X \in \ssf{R}$ for all $m \in \Z$ by definition, and hence $\Gamma_{(m,v)}(X) \cong \Hom_{\Lambda}(P_v,\sigma^{-m}X) \cong (\sigma^{-m}X)\varepsilon_v$ is a finite dimensional right $K_v$-module for all $(m,v) \in \Z Q$. Then $C$ must be locally finite as well. As $\E$ is projective (P), $C$ is also upper bounded. If the upper bound is $(d,v_d) \in \Z Q$, then $C = \bigoplus_{(m,v) \in J} C_{(m,v)}$, where $J = \{ (m,v) \in \Z Q \mid (m,v) > (n,w),~ (m,v) \not> (d,v_d) \}$. The set $J$ is finite as $J=((n,w),(d,v_d)] \cup \{(m,v) \in \Z Q \mid (m,v)>(n,w),~(m,v)\vert\vert (d,v_d)\}$ (c.f. Remark 3.15), so it follows that $C$ is a finite direct sum of finite dimensional modules. Thus $C$ is finitely generated, and hence $\E$ is coherent (C).

    (A): Fix $X \in \H$ and $(m,v) \in \Z Q$. We choose $n > m$ such that $(m,v) \prec (n,w)$ for all $w \in Q_0$. By Lemma 7.2(3) we may choose $n$ large enough such that $\sigma^{-n}X \in \ssf{P} \vee \ssf{R} \subset \mod \Lambda$. We then take an injective resolution $0 \to \sigma^{-n}X \to E \to E'$ of $\sigma^{-n}X$, with $E' = \bigoplus_{l=1}^r E_{v_l}$. This resolution determines a triangle in $D^b(\Lambda)$, which after applying $\sigma^n$ and rotating produces the triangle
    \[\sigma^nE'[-1] \to X \to \sigma^nE \xrightarrow{[1]}.\]
    As $\sigma^nE \in \ssf{I} \subset \CC^{\leq -1}$, this triangle induces a surjection $\sigma^nE'[-1]= \bigoplus_{l=1}^r\sigma^{n+1}P_{v_l} \to X$ in $\H$. We have $(m,v) \prec (n+1,v_l)$ for $1 \leq l \leq r$ by construction, so $\E$ is ample (A).
\end{proof}

Theorem 7.3(3) tells us that $\PP(\cS)$ should be thought of as a noncommutative projective curve, and we view Theorem 7.3(2) as a version of Beilinson's derived equivalence \cite{beilinson_coherent_1978} for these curves. 

The \textit{preprojective algebra} $\Pi(\cS)=\Pi$, introduced for bimodule species by Dlab and Ringel in \cite{dlab_preprojective_1980}, plays an important role in the representation theory of the species $\cS$. The \textbf{preprojective algebra} may be defined as the orbit algebra $\Pi = \bigoplus_{n \geq 0} \Hom_{\Lambda}(\Lambda,\tau^{-n}\Lambda)$ of $\Lambda$ under $\tau^-$, where we write $\tau^{-n}$ to mean $(\tau^-)^n$. The multiplication on $\Pi$ is given by $r \times_{\Pi} s := \tau^{-n}(r) s$ for $r \in \Pi_m$ and $s \in \Pi_n$. As $\sigma^{-1}(\ssf{P}) \subset \ssf{P}$, $\Pi$ is isomorphic to the orbit algebra $\bigoplus_{n \geq 0} \Hom_{\H}(\Lambda,\sigma^{-n}\Lambda)$.

Before we can relate $A$ to the preprojective algebra $\Pi$, we need to recall the notion of \textit{quasi-Veronese} algebras from \cite{mori_construction_2013}. If $I$ is a set with a partition $I = \bigsqcup_{j \in J} I_j$ such that each $I_j$ is finite, and if $B = \bigoplus_{i,i' \in I} B_{ii'}$ is an $I$-indexed algebra, the \textbf{quasi-Veronese} $B^{[J]}$ is the $J$-indexed algebra $B^{[J]} = \bigoplus_{j_0,j_1 \in J} (\bigoplus_{i_0 \in I_{j_0},i_1 \in I_{j_1}}B_{i_0i_1})$ with multiplication 
\[(a_{i_0i_1}) \times (b_{i_1'i_2}) = \bigg(\sum_{i_1=i_1'} a_{i_0i_1}b_{i_1'i_2}\bigg) \in B^{[J]}_{j_0,j_2}\]
for $(a_{i_0i_1}) \in B^{[J]}_{j_0,j_1}$ and $(b_{i_1'i_2}) \in B^{[J]}_{j_1,j_2}$. The maps $\bigoplus_{i \in I} M_i \mapsto \bigoplus_{j \in J}(\bigoplus_{i \in I_j} M_{i})$ determine an exact equivalence $\Gr B \simeq \Gr B^{[J]}$. There is a natural partition 
\[ \Z \times Q_0=\bigsqcup_{m \in \Z}\{(m,v)\mid v \in Q_0\},\] 
so we may construct the quasi-Veronese $A^{[\Z]} = \bigoplus_{m,n \in \Z} (\bigoplus_{v,w \in Q_0} A_{(m,v)(n,w)})$. Note that the positive $\Z Q$-indexing on $A$ implies $A^{[\Z]}$ is positively $\Z$-indexed. 

\begin{lemma}
    There is an isomorphism of $\Z$-indexed algebras $A^{[\Z]} \cong \mathbf{B}_{\Z}\Pi$.
\end{lemma}
\begin{proof}
    It follows from (7.1) that we have isomorphisms $\sigma:A^{[\Z]}_{mn} \cong A^{[\Z]}_{m+1,n+1}$ for all $m,n \in \Z$. Then $A^{[\Z]} \cong \mathbf{B}_{\Z}S$ by \cite[Lemma 2.4]{van_den_bergh_noncommutative_2011}, where $S= \bigoplus_{n \in \Z} A^{[\Z]}_{0n}$ is a $\Z$-graded ring with multiplication $s_m \times_S s_n := s_m \sigma^m(s_n) \in S_{m+n}$ for $s_m \in S_m,~s_n\in S_n$. Then, as $S_n = \Hom_{\H}(\sigma^n\Lambda,\Lambda) \cong \Hom_{\H}(\Lambda,\sigma^{-n}\Lambda)$, we have an isomorphism of graded abelian groups $\varphi:S \xrightarrow{\sim} \Pi$ sending $s_n \in S_n$ to $\sigma^{-n}(s_n) \in \Pi_n$. This isomorphism is compatible with the multiplications, as
    \begin{align*}
        \varphi(s_m \times_S s_n) = \sigma^{-(m+n)}(s_m\sigma^m(s_n)) &= \sigma^{-n}(\sigma^{-m}(s_m))\sigma^{-n}(s_n)\\
        &=\sigma^{-n}(\varphi(s_m))\varphi(s_n)\\
        &= \varphi(s_m) \times_{\Pi}\varphi(s_n).
    \end{align*}
    Thus $S$ and $\Pi$ are isomorphic as $\Z$-graded rings, and hence $A^{[\Z]} \cong \mathbf{B}_{\Z} S \cong \mathbf{B}_{\Z}\Pi$.
\end{proof}

Theorem 7.3(1) and the chain of exact equivalences $\Gr A \simeq \Gr A^{[\Z]} \simeq \Gr \Pi$ ensures that $\Pi$ is coherent and that there is an exact equivalence $\PP(\cS) \simeq \qgr \Pi$. Thus we obtain a generalisation of a known result for algebraic species \cite[Corollary 5.4]{minamoto_ampleness_2012}.

\begin{cor}
    The preprojective algebra $\Pi(\cS)$ of any bimodule species $\cS$ is coherent, and $D^b(\Lambda(\cS)) \simeq D^b(\qgr \Pi(\cS))$.
\end{cor}

We end by showing that $\PP(\cS)$ satisfies a number of properties one may expect from a projective curve, similarly to \cite[Corollary 5.5, Theorem 5.6]{chan_representation_2019}. We let $\ssf{F}$ and $\ssf{T}$ denote the images in $\PP(\cS)$ of the full subcategories $\ssf{I}[-1] \vee \ssf{P}$ and $\ssf{R}$ of $\H$, respectively. We view objects in $\ssf{T}$ as torsion sheaves and objects in $\ssf{F}$ as torsion-free sheaves on $\PP(\cS)$. 

\begin{propn}
    Let $\MM$ be a coherent sheaf on $\PP(\cS)$. \begin{enumerate}
        \item The pair $(\ssf{T},\ssf{F})$ is a torsion pair in $\PP(\cS)$.
        \item $\MM$ splits as a direct sum $\MM \cong \TT_{\MM} \oplus \F_{\MM}$, where $\TT_{\MM} \in \ssf{T}$ and $\F_{\MM} \in \ssf{F}$.
        \item (Grothendieck splitting) $\ssf{F}$ is closed under extensions.
        \item (Serre duality) For all $(m,v) \in \Z Q$ and $l \in \{0,1\}$, there is an isomorphism of $K_v$-vector spaces
        \[ \E xt^l_A(\cP_{(m,v)},\MM)^{\ast} \cong \E xt^{1-l}_A(\MM, \cP_{(m+1,v)}). \]
        \item (Serre finiteness/vanishing) $\E xt^l_A(\cP_{(m,v)},\MM)$ has finite dimension over $K_v$ for $l \geq 0$, and vanishes for $m\gg 0$ if $l>0$.
    \end{enumerate}
\end{propn}
\begin{proof}
    As $\H=\ssf{I}[-1]\vee \ssf{P}\vee \ssf{R}$, it follows from \cite[Lemma 3.3(3)]{chan_representation_2019} and \cite[Theorem 2.1]{dickson_torsion_1966} that $(\ssf{R},\ssf{I}[-1]\vee \ssf{P})$ is a torsion pair in $\H$. Then (1) follows as $\ssf{T}=\QGamma_{\ast}(\ssf{R})$ and $\ssf{F}=\QGamma_{\ast}(\ssf{I}[-1]\vee \ssf{P})$. In light of (1), both (2) and (3) are basic results of torsion theory (see \cite[Section 2]{dickson_torsion_1966}).
    
    To show (4) and (5) it suffices to prove the equivalent statements in $\H$. Both cases for (4) are similar, so we only show for $l=0$. Take $X \in \H$, and without loss of generality we suppose $X$ is indecomposable. As $\sigma$ is an automorphism we have $\Hom_{\H}(\sigma^m P_v,X) \cong \Hom_{\H}(P_v,\sigma^{-m}X)$. Note that $\Hom_{\H}(\sigma^mP_v,X)^{\ast} \cong \Ext^1_{\H}(X,\sigma^{m+1}P_v)$ if $\sigma^{-m}X = E[-1] \in \ssf{I}[-1]$, as we have right $K_v$-linear isomorphisms
    \[ \Hom_{\H}(P_v,\sigma^{-m}X) \cong \Hom_{D^b(\H)}(P_v,\sigma^{-m}X) \cong \Hom_{D^b(\Lambda)}(P_v, \sigma^{-m}X) \cong \Ext_{\Lambda}^{-1}(P_v,E) = 0, \] 
    and similarly we have left $K_v$-linear isomorphisms
    \[ \Ext^1_{\H}(X,\sigma^{m+1}P_v) \cong \Hom_{D^b(\H)}(\sigma^{-m}X,\sigma P_v[1]) \cong \Hom_{D^b(\Lambda)}(E[-1],E_v) \cong \Ext_{\Lambda}^1(E,E_v) = 0\] 
    as $E_v$ is an injective $\Lambda$-module. We now suppose $\sigma^{-m}X \in \ssf{P} \vee \ssf{R} \subset \mod \Lambda$, in which case we have right $K_v$-linear isomorphisms
    \[ \Hom_{\H}(P_v,\sigma^{-m}X) \cong \Hom_{D^b(\H)}(P_v,\sigma^{-m}X) \cong \Hom_{D^b(\Lambda)}(P_v,\sigma^{-m}X) \cong \Hom_{\Lambda}(P_v,\sigma^{-m}X). \]
    Then (4) follows from the following sequence of left $K_v$-linear isomorphisms:
    \begin{align*}
        \Hom_{\H}(\sigma^m P_v,X)^{\ast} &\cong \Hom_{K_v}(\Hom_{\Lambda}(P_v,\sigma^{-m}X),K_v)\\
        &\cong \Hom_{K_v}(\sigma^{-m}X\otimes_{\Lambda} \Lambda\varepsilon_v, K_v)\\
        &\cong \Hom_{\Lambda}(\sigma^{-m}X,\Hom_{K_v}(\Lambda \varepsilon_v,K_v))\\
        &= \Hom_{\Lambda}(\sigma^{-m}X,\varepsilon_vD\Lambda) \\
        &\cong \Hom_{D^b(\Lambda)}(\sigma^{-m}X,\sigma P_v[1]) \\
        &\cong \Hom_{D^b(\H)}(X,\sigma^{m+1}P_v[1]) \\
        &\cong \Ext^1_{\H}(X,\sigma^{m+1}P_v).
    \end{align*}
    
    As the sequence $\E = (\sigma^mP_v)_{(m,v) \in \Z Q}$ is ample, any $X \in \H$ has a resolution by sums of the $\sigma^mP_v$, and so to prove local finiteness for (5) it suffices to show for $X=\sigma^mP_v$. We also only need to consider $l=0,1$ by Theorem 7.3(3). When $l=0$ we have $\Hom_{\H}(\sigma^nP_w,\sigma^mP_v) = A_{(m,v)(n,w)}$, and so local finiteness follows from $(\ast)$ and Lemma 3.8. When $l=1$, local finiteness then follows from Serre duality. In fact, Serre duality shows 
    \[ \Ext_{\H}^1(\sigma^mP_v,X)^{\ast} \cong \Hom_{\H}(X,\sigma^{m+1}P_v) \cong \Hom_{\H}(\sigma^{-(m+1)}X,P_v)\]
    for every $X \in \H$. Then $\Ext^1_{\H}(\sigma^mP_v,X)$ vanishes for $m \gg 0$, as $\sigma^{-(m+1)}X \in \mod \Lambda$ for $m \gg 0$ by Lemma 7.2(3) and there are only finitely many $N$ in $\mod \Lambda$ with $\Hom_{\Lambda}(N,P_v) \neq 0$.
\end{proof}

In the case that $\cS$ has underlying quiver $Q=\widetilde{A}_1$, many well-known results about $\qgr\Pi(\cS)$ for algebraic species were recovered in \cite{chan_species_2016}. For example, the algebra $A$ is noetherian of global dimension 2 \cite[Theorem 5.2, Proposition 6.1]{chan_species_2016}, generalising a result for the algebraic case \cite[Proposition 4.2]{baer_tame_1987}, and a classification of the regular $\Lambda$-modules is given \cite[Theorem 10.17]{chan_species_2016}, mirroring a result of Ringel requiring $\mathrm{rad}\Lambda$ to be non-simple \cite[Theorems 4]{ringel_representations_1976}. Similarly, studying $\PP(\cS)$ when $\cS$ is \textit{tame} ($Q$ is extended Dynkin) may allow one to extend known results for tame algebraic species to the non-algebraic setting, providing new insights into the representations of non-algebraic tame species.

\section*{Acknowledgements}
We wish to thank Daniel Chan for his many helpful comments and suggestions.




\end{document}